\documentclass[10pt]{article}
\usepackage[latin1]{inputenc}
\usepackage{graphicx}
\usepackage{amssymb}
\usepackage{amsfonts}
\usepackage[thmmarks]{ntheorem}
\usepackage[all]{xy}

\def\RR{{\mathbb R}}
\def\CC{{\mathbb C}}

\def\bra{\langle}
\def\ket{\rangle}
\def\id{{\rm Id}}
\def\tr{\mbox{\rm Tr}}
\def\ker{\mbox{\rm Ker}}
\def\im{\mbox{\rm Im}}

\def\bea{\begin{eqnarray}}
\def\eea{\end{eqnarray}}

\def\be{\begin{equation}}
\def\ee{\end{equation}}

\newtheorem{theorem}{Theorem}
\newtheorem{definition}{Definition}
\newtheorem{lemma}{Lemma}
\newtheorem{cor}{Corollary}
\newtheorem{propo}{Proposition}
\theoremstyle{nonumberplain}
\theorembodyfont{\normalfont}
\theoremseparator{:}
\theoremsymbol{$\P$}
\newtheorem{demo}{Proof}
\makeatletter
\newcommand{\manuallabel}[2]{\def\@currentlabel{#2}\label{#1}}
\makeatother

\newenvironment{remark}{\small {\bf Remark:}}{}

\manuallabel{sec54}{3.4}

\begin{document}
\title{On the definition of spacetimes in Noncommutative Geometry: Part I}
\author{Fabien Besnard, Nadir Bizi}
\maketitle

\begin{abstract}
In this two-part paper we propose an extension of Connes' notion of even spectral triple to the Lorentzian setting. This extension, which we call a \emph{spectral spacetime}, is discussed in part II where  several natural examples are given which are not covered by the previous approaches to the problem. Part I only deals with the commutative and continuous case of a manifold.  It contains all the necessary material for the generalization to come in part II, namely the characterization of the signature of the metric in terms of a time-orientation 1-form and a natural Krein product on spinor fields. It turns out that all the data available in Noncommutative Geometry (the algebra of functions, the Krein space of spinor fields, the representation of the algebra on it, the Dirac operator, charge conjugation and chirality), but nothing more, play a role in this characterization.  Thus, only space and time oriented spin manifolds of even dimension are considered for a noncommutative generalization in this approach. We observe that these are precisely the kind of manifolds on which the modern theories of spacetime and matter are defined.
\end{abstract}

\tableofcontents
\section{Introduction}
Noncommutative geometry, as initiated by Alain Connes, is an operator algebraic framework which generalizes Riemannian manifolds in a way which harmoniously gathers continuous and discrete spaces, as well as truly noncommutative examples. This generalization happens to be just enough to allow the recasting of the Standard Model of particle physics as a noncommutative Kaluza-Klein theory. This last application, thus far limited to the compact Riemannian case, has been a strong motivation for the search of a semi-Riemannian  extension of Connes' noncommutative geometry.

There  have been several important attempts in that direction already, but this endeavour is plagued with   difficulties. Some are analytical: among them we find the noncompactness issue and the definition of the spectral action. This last problem is probably the most elusive for the time being. A different kind of problem is the characterization of the physical Lorentzian signature among the general semi-Riemannian ones. This has been investigated in particular in \cite{bes1}, \cite{franco}. The present paper exclusively deals with this last issue. We will present a solution which, as we will see, is   different from existing ones on some essential points. However we will argue that this solution   encompasses both discrete and continuous spaces, as it should, whereas its predecessors did not.

However, it should be stressed that the solution we will propose is incomplete. What we will do is discuss the commutative and continuous case, and use it to motivate the general definition of a noncommutative structure which should correspond to Lorentzian spacetimes in the same way spectral triples correspond to Riemannian manifolds. But in this definition we will totally  ignore the analytical aspects. Hence, in full rigor, it is only applicable to the finite-dimensional (hence discrete) case. There  are two reasons why we are being so careless. The first and most important is to keep this work within reasonable bounds. The second is that we think it is more urgent to present the reader with a careful motivation of the general idea, followed by several examples, in order to let him/her ponder the relevance of our approach,  than elaborate  the subtler aspects of the question. These would only set us out of focus, as well as, maybe, repell readers with a physicist background, who generally tend to shun discussions which they think are too formal\footnote{Actually, we are also mainly interested in the potential physical  applications of noncommutative geometry, and  we believe that in the quantum  gravitational regime  the analytical questions  will essentially go away.}. That said, we believe that it should be possible to use the work \cite{vddpr} by van den Dungen, Paschke and Rennie to complete the work presented here. In fact our results are incomplete in yet another way: we only deal here with the even (KO-)dimensional case.

Let us now be more specific about what we are going to do. We will argue that in order to describe Lorentzian manifolds, we must replace a spectral triple  $({\cal A},{\cal H},D)$ with a structure $({\cal A},{\cal K},D)$ which we call a \emph{spectral spacetime}. Though we will not describe in this introduction the details of the definition, let us just say first that the Hilbert space ${\cal H}$ appearing in the data of the spectral triple, must be replaced with a Krein space ${\cal K}$. There is no novelty here. We will point out that this requires us to restrict to the class of space-oriented manifolds. Since we are also given a total (space+time) orientation, the manifolds which can be described in a spectral way are space and time oriented ones. More importantly, we will propose a definition of  time orientation in the noncommutative case: such an orientation will be defined to be a noncommutative 1-form $\beta$ which must be real in some sense and turn the Krein product into a scalar product. Much more importantly perhaps, we will argue that we must \emph{not} ask ${\cal A}$ to be a $C^*$-algebra, in fact not even a $*$-algebra, just an algebra of operators acting on ${\cal K}$. This step is not as radical as one might think: it is the natural consequence of the downplay of the Hilbert structure in front of the Krein structure.  However, without a $C^*$-structure we take the risk of not being able to reconstruct a space out of the algebra, and this is why we will propose a natural axiom of \emph{reconstructibility}, which is that the algebra is closed under the Hilbert adjunction defined by at least one orientation form $\beta$. In such a case the algebra ${\cal A}$ may posses many different, yet isomorphic, $C^*$-structures: one for each orientation form satisfying the reconstructibility condition.

As we have alluded to above, the main arguments we will give in favour of our notion of spectral spacetime are the continuous case on one hand, and discrete examples on the other. These discrete examples are constructed over a finite graph. The first, which we call the \emph{canonical antilorentzian spectral spacetime over a positively weighted graph} is just an appropriate Wick rotation of the usual spectral triple one can build on such a graph to recover the geodesic distance thanks to Connes' distance formula. It is easily seen to be  a spectral spacetime but none of the structures proposed before by various authors. The second example is an elaboration of the first, which we decorate with noncommutative algebras on the vertices, and discrete parallel transport operators on the edges. We call it the  \emph{split Dirac structure} over a graph.  The Dirac operator in this structure turns out to have an interesting relation with the discretization of the Dirac operator proposed by Marcolli and van Suijlekom in \cite{MVS}: the latter is obtained   by acting with the former on special ``graph states'' and projecting out the result on the space of these states. Moreover, we will see that for  the Dirac structure to be a spectral spacetime it is necessary that the discrete parallel transport operators satisfy relations which are exactly\footnote{There is actually a small difference, which does not show up in the most natural cases, and is not worth explaining in this introduction.} those  required of a spin connection on a manifold. Finally, this example will also show us that the reconstructibility condition, however natural it might be, is quite restrictive in the noncommutative case since it is equivalent to the existence of a parallel (i.e. covariantly constant)  timelike vector field.

We have approached our subject in a very slow and (we hope) careful way keeping these general principles in mind while writing the text: stack definitions, lemmas and propositions until the time is ripe to state a theorem, stay as geometric as possible, stay as general as possible unless it spoils clarity, always be clear on what convention is used, be as self-contained as possible, display (almost) every computation. The resulting pace is admittedly very slow, and it would certainly have been possible to take some  short cuts.  This is why we will give below a little guide for readers who wish to go directly to the main points. We hope nonetheless that the material which may  be  skipped at first could prove useful on a second read for those who wish to deepen some side question.

The two parts of the paper can essentially be read independently. In the first part, which is devoted to the classical (i.e. continuous and commutative) case, we do not claim much originality. It can be seen   mainly as a motivation for the definition of spectral spacetimes coming in the second part. However, we have tried to give a unified and, we think, new presentation of several subjects (real structures on Clifford algebras, Wick rotations,  canonical hermitian forms on Clifford algebras and  on spinors) which are otherwise treated in isolation from one another in the literature. Moreover we have proven some results which are probably known as ``folk theorems'' by some, but that we could not locate in the literature. Finally the first part contains some technical lemmas used in the second part. The second part provides the definition of spectral spacetimes as well as the study of several   examples, as explained above.

The first part of the paper is organized as follows:   after this introduction, we  make a somewhat long excursion into Clifford algebras in section 2. The goal there will be to give a characterization of the Lorentzian signature of the metric with tools available  in noncommutative geometry. To this end, we first give a   description of Wick rotations in terms of commuting real structures, then we recall the existence of Krein structures on spinor modules which are compatible with a given real structure, and use them to characterize the Wick rotations to Euclidean signature (\emph{Robinson's alternative}, proposition \ref{propo7}). Finally our goal is achieved in theorem \ref{maintheorem}, which characterizes the Lorentzian and anti-Lorentzian signatures by the possibility of using a single vector to turn the Krein structure on a given spinor module into a Hilbert structure. This section also makes the connection between the real structures on the Clifford algebra and the charge conjugation operator which is used in noncommutative geometry, and also contains some extra material on Clifford algebras, e.g. the  $\sigma$-product, which extends a given metric to the whole algebra.

In section 3 we globalize the notions introduced in section 2 to semi-rieman\-nian manifolds: real structures on Clifford bundles, Krein structures on spinor bundle. We point out that Krein structures compatible with a given real structure exist if and only if the manifold is space or time orientable, according to the signature (theorem \ref{theorem2}). We also give a detailed explanation of the way in which the geometric properties of the spin connection are encoded in the algebraic properties of the Dirac operator and charge conjugation.  Finally we recall how the Dirac operator can be Wick rotated, and we translate our local characterization of  Lorentzian signature from the previous section into global terms (theorem \ref{carac2}), namely in terms of an orientation 1-form. We point out that the exactness of this form amounts to a causality condition which is called \emph{stable causality}. This closes part I.

We advise the hurried reader to use the following strategy: read quickly subsection \ref{sec24} for the definition of Krein product on spinor space compatible with a given real structure, then read subsection \ref{sec27} to have theorem \ref{maintheorem} in mind, and skip to section 3.

In the whole paper some notions are defined in the bulk of the text. When it happens, the name of the notion is always italicized.


\section{Local constructions}
\subsection{General definitions and conventions}


In this section we will be interested in the complex Clifford algebra generated by a real vector space $V$ of even dimension $n=2k$ equipped with a non-degenerate quadratic form $Q$ of signature $(p,q)$. The bilinear form associated with $Q$ will be denoted by $B$. Our notations and conventions are summarized below:

\begin{itemize}
\item As is traditional in physics, what we call the signature of a real bilinear symmetric (or complex sesquilinear hermitian) form  is the triple $(n_+,n_-,n_0)$, where $n_0$ is the nullity, $n_-$ the negative and $n_+$ the positive indices of inertia. When the form is non-degenerate, we leave out $n_0=0$ and write $(n_+,n_-)$ instead.
\item We give names to the following cases :  \emph{Lorentz signature} is $(n-1,1)$,   \emph{anti-Lorentz signature} is $(1,n-1)$, \emph{Euclidean signature} is $(n,0)$ and \emph{neutral signature} is $(n,n)$.
\item The real Clifford algebra $Cl(V,Q)$ is  defined using the following convention\footnote{We suggest the name ``anti-Clifford'' for the algebra defined  using the convention $v^2=-Q(v)$. We think that making the sign conventions explicit by using the words anti-Clifford and antilorentzian would be greatly beneficial to the mathematical physics community.}  
$$
v^2=+Q(v)
$$
for all $v\in V$. In later sections, in the context of manifolds, we will put the emphasis on the bilinear form and write $Cl(V,B)$ instead of $Cl(V,Q)$. We let $\CC l(V)=Cl(V,Q)\otimes \CC$ be its complexification, and $V^\CC$ be the complexification of $V$. We consider that $V\subset Cl(V,Q)\subset \CC l(V)$ in the natural way, i.e. we do not write down explicitly the   embedding $i : V\rightarrow Cl(V,Q)$. We still denote by $Q$ the natural extension of $Q$ such that $v^2=Q(v)$ holds for all $v\in V^\CC$.
\item We let $c : \CC l(V)\rightarrow \CC l(V)$ be the real structure defined by $c(a\otimes \lambda)=a\otimes \bar \lambda$.
\item We let $T : \CC l(V)\rightarrow \CC l(V)$ be the unique linear antiautomorphism which restricts to the identity on $V$. It is called {\it the principal anti-involution}. We often write $a^T$ instead of $T(a)$. 
\item We note that $c\circ T=T\circ c$ is the unique antilinear antiautomorphism of $\CC l(V)$ which extends $\id_V$. We will   write $a^\times=c(a^T)$.
\item We let $\gamma$ be the {\it principal involution}, that is, the unique automorphism of $\CC l(Q)$ which extends $-\id_V$. This a grading operator which decomposes the Clifford algebra into the sum $\CC l(V)=\CC l^0(V)\oplus \CC l^1(V)$ of its even and odd parts. 
\item For any $a$ and invertible $g$ we write $L_g(a)=ga$, $R_g(a)=ag$, $Ad_g(a)=gag^{-1}$. The (complex) Clifford group is defined by $\Gamma_\CC=\{g\in\CC l(V)|$ $Ad_g(V^\CC)\subset V^\CC\}$, the Pin group is $Pin(Q)=\{g\in \Gamma_\CC|c(g)=g$ and $gg^T=\pm 1\}$, and the Spin group is $Spin(Q)=Pin(Q)\cap Cl(V,Q)^0$. Remember that the Clifford group is generated by non-isotropic vectors, and that its elements satisfy $gg^\times\in\RR$.
\item Let $e_1,\ldots,e_{n}$ be a pseudo-orthonormal basis of $V$. Then we denote by $\omega=e_1\ldots e_n$ the volume element. It depends on the pseudo-orthonormal basis chosen only up to a sign, which we can fix by chosing an orientation of $V$. The volume element anticommutes with every odd element of $\CC l(V)$. It has the following properties :
$$\omega^T=(-1)^{{n\over 2}}\omega, \omega^2=(-1)^{{n\over 2}+q}$$
 
\item If $\rho$ is a spinor representation, we set  $\chi:=(-i)^{{n\over 2}+q}\rho(\omega)$. It is called the chirality operator and  always satisfies $\chi^2=1$. When the representation is $c$-admissible (to be defined later), it will also satisfy $\chi^\times=(-1)^q\chi$. 
\end{itemize}

Our general reference on Clifford algebras is \cite{crum}. For future use we note the following fact : let $\phi : \CC l(V)\rightarrow \CC l(V)$ be an automorphism or antiautomorphism which stabilizes $V^\CC$. Then $\phi\circ T$ and $T\circ \phi$ coincide on $V^\CC$, hence on the whole algebra.  Thus $T$ and $\phi$ commute. 



\subsection{Real structures and local Wick rotations}

What we wish to do in this section is to understand Wick rotations algebraically. We are given at the start  the vector space $V$ equipped with a quadratic form $Q$. The piece of data $(V,Q)$ is equivalent to $Cl(V,Q)$. We stress that the Clifford algebra is not to be seen only as a real algebra, but as a real algebra equipped with a particular set of generators, namely $V$. Since we want to change the signature of the quadratic form, it is natural to embed $Cl(V,Q)$ in its complexification $\CC l(V)$ which will remain constant when $Q$ is varied. Recovering $Q$ from $\CC l(V)$ amounts to fix a particular real form for the complex algebra $\CC l(V)$, by  way of a \emph{real structure}, i.e.  an involutive antilinear automorphism which stabilizes $V^\CC$.  The real structure which has $Cl(V,Q)$ as its set of fixed points will be denoted by $c$ throughout the text. Clearly, the data $(V,Q)$ and $(\CC l(V),c)$ are equivalent. Given a general real structure $\sigma$ we define:

%

\begin{itemize}
\item The real subspace of $\sigma$-real vectors $V_\sigma:=\{v\in V^\CC|\sigma(v)=v\}$.
\item The map $u_\sigma : v\mapsto {v+\sigma(v)\over 2}+i{v-\sigma(v)\over 2}$, which is easily seen to be an isomorphism of real vector spaces from $V$ onto $V_\sigma$.
\item The bilinear form $B_\sigma(v,v'):={1\over 2}(B(\sigma(v),v')+B(\sigma(v'),v))$ on $V$, and its associated quadratic form $Q_\sigma(v):=B(\sigma(v),v)$.
\end{itemize}

Since $\sigma(v)$ does not in general belong to $V$ it is not obvious at first sight  that $B_\sigma(v,v')$ is real. However we can observe that $B_\sigma(v,v')={1\over 2}(\sigma(v)v'+v'\sigma(v)+\sigma(v')v+v\sigma(v'))$ and on this form it is clear that $\sigma(B_\sigma(v,v'))=B_\sigma(v,v')$, hence this complex number is in fact real. 

We also note that the restriction of the quadratic form $Q$ to $V_\sigma$ is real, since   $\sigma(w^2)=\sigma(w)^2=w^2$ for every $w\in V_\sigma$. Furthermore, since $V_\sigma+iV_\sigma=V^\CC$, the complex algebra generated by $V_\sigma$ is $\CC l(V)$, which we can then identify with $Cl(V_\sigma,Q_{|V_\sigma})\otimes \CC$.  In particular we note  that $Q_{|V_\sigma}$ is non-degenerate.

Finally the calculation $Q(u_\sigma(v))=u_\sigma(v)^2={1\over 2}(v\sigma(v)+\sigma(v)v)=Q_\sigma(v)$ shows that $u_\sigma$ is an isometry from $(V,Q_\sigma)$ to $(V_\sigma,Q)$.

Hence we have defined from $\sigma$ a non-degenerate quadratic form $Q_\sigma$ on $V$ and a real subspace $V_\sigma$ of $V^\CC$ such that $V_\sigma\oplus iV_\sigma=V^\CC$ and $Q_{|V_\sigma}$ is real.


This construction can be inversed. Consider a $n$-dimensional  real vector space $W$ of $V^\CC$ such that $Q$ is real on $W$. Then   $W\oplus iW=V^\CC$ (indeed, if there exists a nonzero $w\in W$ such that $iw\in W$ then $Q(w+iw)$ is not real). Define $\sigma : V^\CC\rightarrow V^\CC$ to be antilinear and satisfy $\sigma_{|W}=\id_W$. It is easy to check that $\sigma(v_1v_2+v_2v_1)=\sigma(v_1)\sigma(v_2)+\sigma(v_2)\sigma(v_1)$, hence $\sigma$  extends as an antilinear algebra automorphism. 






Hence we have the following one-to-one correspondence:

\bea
\{\mbox{ real structures on }\CC l(V)\}&\simeq &\{n-\mbox{dimensional real subspaces }\cr
&&W\subset V^\CC\mbox{ such that }Q_{|W}\mbox{ is real}\}\nonumber\cr
\eea




Since $c$ is a fixed ``background'' structure, it is natural to be particularly interested in real structures which commute with it, and we will call them {\it admissible real structures}. They turn out to correspond to   Wick rotations of the quadratic form $Q$.

\begin{lemma}\label{decadix} The following are equivalent. 
\begin{enumerate}
\item\label{d3} The  real structure $\sigma$ is admissible.
\item\label{d4} The subspace $V_\sigma$ is stable by $c$.
\item\label{d5} The subspace $V$ is stable by $\sigma$.
\item\label{d2} The real structure $\sigma$ restricts to a $Q$-orthogonal symmetry of $V$.
\item\label{d1} The subspaces $V_+:=V\cap V_\sigma$ and $V_-:=V\cap iV_\sigma$ form a  $Q$-orthogonal decomposition of $V$.
\end{enumerate}
\end{lemma}
\begin{demo}
It is immediate that (\ref{d3})$\Rightarrow $(\ref{d5}), and for the converse it suffices to observe that $\sigma\circ c$ and $c\circ \sigma$ are two algebra automorphisms of $\CC l(V)$ which coincide on $V$. The equivalence between (\ref{d3}) and (\ref{d4}) is obtained by symmetry. 

Of course (\ref{d2})$\Rightarrow $(\ref{d5}) is trivial, and  to see that (\ref{d5})$\Rightarrow $(\ref{d2}) we observe that $Q(\sigma(v))=\sigma(v)^2=\sigma(v^2)=\sigma(Q(v))=Q(v)$, hence $\sigma$ restricts to a $Q$-orthogonal transformation of $V$ which is moreover involutive.

Finally we see that $V_+$ is the $+1$-eigenspace and $V_-$ is the $-1$-eigenspace of the $\RR$-linear operator $\sigma_{|V}$, hence (\ref{d1}) and (\ref{d2}) are equivalent.
\end{demo}

Let us consider an admissible real structure $\sigma$. Using the decomposition $V=V_+\oplus V_-$ given by point \ref{d1} above, we see that if $v=v_++v_-$, with $v_\pm\in V_\pm$, then the isometry
$$u_\sigma : (V,Q_\sigma)\longrightarrow (V_\sigma,Q_{|V_\sigma})$$
takes the simple form

$$v\longmapsto v_++iv_-$$

Hence $u_\sigma$ ``puts an $i$'' in front of the elements of $V_-$: this is what is called a \emph{Wick rotation}.  Clearly $Q_\sigma$ is positive definite iff  $Q$ is positive definite on $V_+$ and negative definite on $V_-$. In this latter case we say that $u$ is a rotation ``to Euclidean signature''.

As the proposition below shows, admissible real structures can be expressed in terms of particular elements of the Clifford group. 

\begin{propo}\label{admreal}
\begin{enumerate} 
\item The  real structures are of the form $\sigma=Ad_b\circ c$ with $b\in \Gamma_\CC$ such that $bc(b)=\lambda\in\RR$. We also have  $b^T=\alpha b$, $\alpha=\pm 1$, $b^\times=\alpha c(b)$.
\item The real structure  $\sigma=Ad_b\circ c$  with $b\in \Gamma_\CC$ is admissible iff $c(b)=e^{i\theta}b$, with $\theta\in\RR$. In this case we can choose $b$ to satisfy $c(b)=b$ and  $b^2=\lambda=\pm 1$, in which case we say that it is  \emph{real and normalized}. Then  $b^\times=\alpha b$, $b$  belongs to $Pin(Q)$, and is unique up to a sign.
\end{enumerate}
\end{propo}
\begin{demo}
It is obvious that  $\sigma$ is an antilinear automorphism which preserves $V^\CC$  when it is of the form $Ad_b\circ c$ with $b$ in the Clifford group. If $bc(b)=\lambda\in\RR$ it is moreover an involution since $\sigma^2(a)=bc(b)ac(b)^{-1}b^{-1}=a$ for all $a\in \CC l(V)$.

Conversely, if $\sigma$ is a real structure, then $\sigma\circ c$ is an automorphism of the Clifford algebra which preserves $V^\CC$, and it is then of the form $Ad_b$ with $b\in \Gamma_\CC$. Hence $\sigma=Ad_b\circ c$. Since $\sigma$ is an involution we have  $bc(b)=\lambda\in\CC$ by the calculation above. Since $c(b)$ is then equal to $b^{-1}$ up to a constant, it commutes with $b$, from which we obtain that $\lambda$ is real. 

 The other properties follow from the ones just proved: since $b\in\Gamma_\CC$, $bb^\times$ is a constant, hence $c(b)b^T$ is a constant, and from $c(b)b=\lambda\in \RR$ we get that $b$ is proportional to $b^T$. The involutory property of $T$ forces the proportionality constant to be a sign $\alpha$. Then $b^\times=c(b^T)=\alpha c(b)$.

Now it is easy to check that $\sigma$ commtutes with $c$  iff  $c(b)=e^{i\theta}b$ for some $\theta\in \RR$. Then  $b'={e^{i\theta/2}\over\sqrt{|\lambda|}}b$ is normalized and real and one has $Ad_b=Ad_{b'}$. If $b''$ is another normalized and real element such that $Ad_{b''}=Ad_{b}$ then $b''=\mu b'$ with $\mu\in\CC$. From reality one has $\mu\in\RR$, and from normalization one has $\mu^2=1$.

Since $b'$   satisfies $b'(b')^\times =\pm 1$ and is real it is in the Pin group.
\end{demo}


{\small At this point one might like to have an example of a non-admissible real structure. For this, consider $\RR^2$ with an Euclidean metric and $(e_1,e_2)$ an orthonormal basis. Then for $t\in\RR$, let $b_t$ be the Clifford group element $b_t=\cosh t+i(\sinh t)e_1e_2$. Then $c(b_t)=b_{-t}=b_t^{-1}$. Hence $b_t$ satisfies the hypotheses of the first part of proposition \ref{admreal}, but not the second (except if $t=0$). If we denote by $V$ the vector of components of $v\in V^\CC$ in the chosen basis, then the real structure $\sigma_t=Ad_{b_t}\circ c$ is given matricially by $V\mapsto O_t\bar V$ where $O_t$ is the matrix

$$O_t=\pmatrix{\cosh 2t& i \sinh 2t\cr -i\sinh 2t& \cosh 2t }$$ 

One then easily shows that the metric $B_{\sigma_t}$ associated to the real structure $\sigma_t=Ad_{b_t}\circ c$ is $\cosh 2t$ times the original metric. More generally, one can show in (even) dimension $n$ and in the Euclidean case that there always exists an orthonormal basis of $V$ in which a general real structure is given by a matrix which is a direct sum of $I_p$, $-I_q$, and $2\times 2$ blocks $\pm O_{t_k}$ with $O_{t_k}$ as above.
}

\subsection{Hermitian forms on the complex Clifford algebra}\label{sec23}

Wick rotations to Euclidean signature are very special. This can be best seen by extending the quadratic form $Q_\sigma$ to the whole Clifford algebra. The automorphism $\sigma\circ T$ will then appear naturally as an adjunction, and its exceptional character when $Q_\sigma$ is positive definite will be made manifest.

Let $(e_i)_{1\le i\le n}$ be a pseudo-orthonormal basis of $V^\CC$. For any subset $I\subset \{1;\ldots,n\}$ we write $e_I=e_{i_1}\ldots e_{i_k}$ where $i_1<\ldots<i_k$ are the elements of $I$. We know that $(e_I)_{I\subset \{1;\ldots,n\}}$ is a basis of $\CC l(V)$. Hence there is a projection map $\tau_n : \CC l(V)\rightarrow \CC$ which sends an element of $\CC l(V)$ to its coordinate on $1=e_\emptyset$. The map $\tau_n$ is called the {\it normalized trace}, a name justified by the following proposition (we refer to \cite{garling} p. 100 for the proof). 

\begin{propo} The projection $\tau_n$ does not depend on the chosen pseudo-ortho\-normal basis. It is the unique linear form $\tau_n : \CC l(V)\rightarrow \RR$ such that $\tau_n(ab)=\tau_n(ba)$ for all $a,b\in \CC l(V)$ and $\tau_n(1)=1$. It also satisfies $\tau_n(a^T)=\tau_n(a)$.
\end{propo}





It is also easy to see that if $a\in\CC l(V)$, then $\tau_n(\gamma(a))=\tau_n(a)$, and if $\sigma$ is a real structure, then $\tau_n(\sigma(a))=\overline{\tau_n(a)}$.
 
\begin{remark}
We can see the above proposition as the reason behind the fact that a product of distinct gamma-matrices is always traceless.
\end{remark}


Given the normalized trace and a real structure, it  is very natural to define

\be
(a,b)_\sigma:=\tau_n(\sigma(a^T)b)\label{KreinCliff}
\ee

for all $a,b\in \CC l(V)$. We call it  \emph{the $\sigma$-product}. It has remarkable properties.

\begin{propo}\label{kreincliff} The $\sigma$-product is  a non-degenerate hermitian form on $\CC l(V)$. The associated quadratic form restricts to $Q_\sigma$ on $V$ and to $Q$ on $V_\sigma$. It satisfies

\be
(w_1\ldots w_k,w_1\ldots w_k)_\sigma=Q(w_1)\ldots Q(w_k)\label{multi}
\ee

for any vectors $w_1,\ldots,w_k\in V_\sigma$. Moreover if $(e_i)_{1\le i\le n}$ is a pseudo-orthonormal basis of $V_\sigma$ for $Q$  then $(e_I)_{I\subset\{1;\ldots;n\}}$ is a pseudo-orthonormal basis of $\CC l(V)$ for the $\sigma$-product.
\end{propo}

\begin{demo}
It is obvious that $(.,.)_\sigma$ is sesquilinear. Moreover we have $(b,a)_\sigma=\tau_n(\sigma(b^T)a)=\tau_n(\sigma(b)^Ta)=\tau_n(a^T\sigma(b))=\overline{\tau_n(\sigma(a^T)b)}=\overline{(a,b)}_\sigma$.

If $v\in V$, then $(v,v)_\sigma=\tau_n(\sigma(v)v)={1\over 2}\tau_n(\sigma(v)v+v\sigma(v))=Q_\sigma(v)$.

Let us prove (\ref{multi}). We have :

\bea
(w_1\ldots w_k,w_1\ldots w_k)_\sigma&=&\tau_n(\sigma(w_k^T)\ldots\sigma(w_1^T)w_1\ldots w_k)\cr
&=&\tau_n(w_k\ldots w_1w_1\ldots w_k)\cr
&=&w_1^2\ldots w_k^2\cr
&=&Q(w_1)\ldots Q(w_k)\nonumber
\eea

Let $(e_i)_{1\le i\le n}$ be a pseudo-orthonormal basis of $V_\sigma$. Then $(e_I,e_I)_\sigma=\pm 1$ by property (\ref{multi}). Moreover if $I\not=J$, $(e_I,e_J)_\sigma=\tau_n(e_I^Te_J)=\pm \tau_n(e_{I\Delta J})=0$. Hence $(e_I)_{I\subset \{1;\ldots;n\}}$ is a pseudo-orthonormal basis of $\CC l(V)$, which shows that $(.,.)_\sigma$ is not degenerate.
\end{demo}


Note that the multiplicative property (\ref{multi}) of the $\sigma$-product does not generalize to elements of $V$. The correct symmetrical statement which apply to elements of $V$ uses the $c$-product instead of the $\sigma$-product.

Here is an elegant (and useful) property of the $\sigma$-product.

\begin{lemma}\label{automatic} Let $\phi : (V,B)\rightarrow (V',B')$ be an isometry between two real vector spaces equipped with nondegenerate bilinear forms. Let $c,c'$ be the canonical real structures on $\CC l(V,B)$ and $\CC l(V',B')$ respectively. Then the isomorphism $\tilde \phi :  \CC l(V,B)\rightarrow \CC l(V',B')$ which canonically extends $\phi$ transforms $(.,.)_c$ into $(.,.)_{c'}$.
\end{lemma}
\begin{demo} For any $a,b\in\CC l(V,B)$ we have
\bea
(\tilde\phi(a),\tilde\phi(b))_{c'}&=&\tau_n'(c'(\tilde\phi(a)^T)\tilde\phi(b))\cr
&=&\tau_n'(c'(\tilde\phi(a^T))\tilde\phi(b)),\mbox{ since }\phi(V)=V'\cr
&=&\tau_n'\circ\tilde\phi (c(a^T)b)\cr
&=&\tau_n(c(a^T)b)=(a,b)_c\nonumber
\eea
where the last step follows from the uniqueness of the normalized trace.
\end{demo}

\begin{remark}
Care must be taken in applying this lemma. For instance if we use it on $u_\sigma$ we obtain  $(\tilde u_\sigma(a),\tilde u_\sigma(b))_\sigma=\tau_n(c(a^T)*_\sigma b)$, where $*_\sigma$ is the Clifford product corresponding to the quadratic form $Q_\sigma$. In particular if $a=b=v$ we obtain $(u_\sigma(v),u_\sigma(v))_\sigma=\tau_n(v*_\sigma v)=Q_\sigma(v)=Q(u_\sigma(v))$ which is correct.
\end{remark}

Taking a $Q$-pseudo-orthonormal basis  $(e_i)_{1\le i\le n}$ in $V$, then the decomposition

\be
\bigoplus_{k=0}^n V^k=Cl(V,Q)\label{decompcliff}
\ee
where  $V^k={\rm Vect}\{e_I|,$ $|I|=k\}$ is immediately seen to be orthogonal for the $c$-product. The following proposition shows in particular that this decomposition is independent of the chosen basis.

\begin{propo}\label{propo4} Let $\phi : (V,B)\rightarrow (V,B)$ be an isometry.  Then $\tilde\phi(V^k)=V^k$ for all $k$.
\end{propo}
\begin{demo}
Clearly $\tilde\phi(V^k)=V^k$ for $k=0,1$. Let us suppose that $\tilde\phi(V^j)=V^j$ for $j\le k$. The sum $V^0\oplus \ldots \oplus V^k\oplus V^{k+1}$ is orthogonal for the $c$-product, hence by the lemma the sum $\tilde\phi(V^0)\oplus\ldots\oplus\tilde\phi(V^k)\oplus \tilde\phi(V^{k+1})=V_0\oplus \ldots\oplus V^k\oplus \tilde\phi(V^{k+1})$ also is. Since we obviously have $\tilde\phi(V^{k+1})\subset V^0\oplus \ldots \oplus V^{k+1}$ we obtain that $\tilde\phi(V^{k+1})=V^{k+1}$.
\end{demo}

Here is another way to understand why the decomposition (\ref{decompcliff}) is independent of the chosen basis. There is a well-known vector space isomorphism $\Theta$ from the exterior algebra $\Lambda V$ to the Clifford algebra which is defined by

\bea
\Lambda V&\longrightarrow&Cl(V,Q)\cr
v_1\wedge\ldots\wedge v_r&\longmapsto&{1\over r!}\sum_\sigma \epsilon(\sigma)v_{\sigma(1)}\ldots v_{\sigma(r)}\label{isotheta}
\eea

Since $\Theta(e_{i_1}\wedge\ldots\wedge e_{i_k})=e_{i_1}\ldots e_{i_k}$ for distinct elements $e_{i_1},\ldots,e_{i_k}$ of the pseudo-orthonormal basis, we see that $V^k=\Theta(\Lambda^k V)$. Moreover, there is a well-known way to extend the bilinear form $B$ to the exterior algebra, which is to decree that $\Lambda^j V$ and $\Lambda^k V$ are orthogonal for $j\not=k$ and to define

$$(u_1\wedge\ldots\wedge u_k,v_1\wedge\ldots\wedge v_k)_B:=\det((B(u_i,v_j))_{1\le i,j\le k})$$

It is left to reader to check that $\Theta$ is an isometry from $(\Lambda V,(.,.)_B)$ to $(Cl(V,Q),(.,.)_c)$. Of course what we have described using $c$   can be extended to a general real structure $\sigma$. However, when passing to the complexification, the reference to the real structure, or quadratic form, vanishes, so that the decomposition

\be
\bigoplus_{k=0}^n V^k\otimes \CC =\CC l(V)\label{decompcliff2}
\ee

is orthogonal for all $\sigma$-products. The vector space $V^k\otimes \CC$ is described in a basis-independent way as $V^k\otimes\CC=\Theta(\Lambda^k V^\CC)$, where we still write $\Theta$ the natural extension of  the isomorphism (\ref{isotheta}) to complex scalars.

\begin{remark}
When $g$ is an element of the Clifford group, $(g,g)_c=\tau_n(g^\times g)=g^\times g$ because $g^\times g$ is real. It is then customary to write $N(g)=g^\times g$ and call it the {\it spinor norm}. 
\end{remark}
\smallbreak

%
%

If $U : \CC l(V)\rightarrow \CC l(V)$ is a linear map, we will write $U^{\times_\sigma}$ for its adjoint relative to the $\sigma$-product. We will also write 
 
$$a^{\times_\sigma}:=\sigma(a^T)$$

for an element $a\in\CC l(V)$. The two notations are consistent thanks to the property 

\be
L_a^{\times_\sigma}=L_{a^{\times_\sigma}}\label{kreinmultadj}
\ee

which is easily checked. We will now see a simple yet important result which tells us how the signature of $(.,.)_\sigma$ depends on that of $Q_\sigma$. We call it \emph{Garling's  alternative}, since the only place where we could locate it is \cite{garling} (p 101), where a direct combinatorial proof is given. We give below a slightly different proof   based on the following   lemma.

\begin{lemma}\label{neutralemma} Let $(K,(.,.))$ be a finite dimensional  space equipped with a non-degenerate hermitian form. If there exists $U\in End(K)$ such that $U^\times U=-\id_K$, then $(.,.)$ is neutral.
\end{lemma}
\begin{demo} We have $(U\psi,U\eta)=(\psi,U^\times U\eta)=-(\psi,\eta)$ for all $\psi,\eta\in K$. Thus if $K=K_+\oplus K_-$ is an orthogonal decomposition of $K$ into subspaces where $(.,.)$ is positive definite and negative definite respectively, we see that $K=UK_+\oplus UK_-$ is an orthogonal decomposition where the signs are swapped. Since $\dim(UK_\pm)=\dim K_\pm$ we conclude by Sylvester's law of inertia.
\end{demo}


\begin{propo}\label{alternative} (Garling's alternative) The $\sigma$-product on $\CC l(V)$ is positive definite whenever $Q_\sigma$ is. It is neutral in every other case.
\end{propo}

\begin{demo} Let $(e_i)_{1\le i\le n}$ be a pseudo-orthonormal basis of $V_\sigma$. Using property (\ref{multi}) of proposition \ref{KreinCliff} on the basis $(e_I)$ we see that $(.,.)_\sigma$ is positive definite if $Q_{|V_\sigma}$ is, and we know that $Q_{|V_\sigma}$ and $Q_\sigma$ have the same signature. 

If $Q_\sigma$ (or equivalently $Q_{|V_\sigma}$) is not positive definite, then there exists $i$ such that $e_i^2=-1$. Let $U=L_{e_i}$. Since $L_{e_i}^{\times_\sigma}=L_{e_i^{\times_\sigma}}=L_{e_i}$, hence $UU^\times=-\id_{\CC l(V)}$ and we conclude by  lemma \ref{neutralemma}.
\end{demo}

We can now characterize Wick rotations to Euclidean signature in terms of the corresponding real structure. The notations are the same as in lemma \ref{decadix}.

\begin{propo}\label{propo6} Let $\sigma$ be an admissible real structure. The following properties are equivalent.
\begin{enumerate}
\item  The couple $(V_+,V_-)$ of $Q$-orthogonal supplementary subspaces of $V$ is such that $Q$ is positive definite on $V_+$ and negative definite on $V_-$.
\item The quadratic form $Q_\sigma$ on $V$ is positive definite.
\item The restriction of the quadratic form $Q$ to $V_\sigma$ is positive definite.
\item  The $\sigma$-product $(.,.)_\sigma$ is positive definite.
\end{enumerate}
\end{propo}



When these properties are satisfied we say that $\sigma$ is an \emph{Euclidean real structure}. Notice that when $\sigma=Ad_b\circ c$ is an Euclidean real structure, and $b$ is normalized, the signs such that $b^\times=\pm b$ and $b^2=\pm 1$ are no longer independent. Indeed, $b^{\times_\sigma}=b^{-1}b^\times b=b^\times$, hence  $bb^\times$   is a positive operator and can only be equal to $1$.  

\begin{remark}
Another specific feature of Euclidean real structures is that the involution $a\mapsto a^{\times_\sigma}$ and the norm $\|a\|_{\infty,\sigma}=\sup_{\|b\|_\sigma=1}\|ab\|_\sigma$ turn  $\CC l(V)$ into a $C^*$-algebra. For more on this, see appendix \ref{cstarstruccliff}.
\end{remark}

\subsection{Krein products on spinor spaces}\label{sec24}

In the rest of this section we let $\rho : \CC l(V)\rightarrow End(K)$ be a  representation of $\CC l(V)$ on a finite dimensional complex vector space $K$ equipped with a non-degenerate hermitian form $(.,.)^K$. We call such a form a  \emph{Krein product}. Remember that  \emph{a  Krein space} (see \cite{bognar}) is a vector space equipped with a a non-degenerate hermitian form $(.,.)$ and at least one operator $\beta$, called a \emph{fundamental symmetry} such that $\beta^2=1$ and $(.,\beta .)$ is positive definite and induces a complete topology. In finite dimension the existence of $\beta$ is obvious, thus $K$ equipped with a Krein product is a Krein space.


We will say that $\rho$ is \emph{$\sigma$-compatible} iff 

$$
\forall a\in \CC l(V),\forall \psi,\phi\in K,  (\rho(a)\psi,\phi)=(\psi,\rho(a^{\times_\sigma})\phi)
$$

This is clearly equivalent to ask  $\rho(v)$  to be Krein-self adjoint for all $v\in V$. When the Krein product is $\sigma$-compatible we will denote by $A^{\times_\sigma}$  the Krein adjoint of $A\in End(K)$.

We can always build a $\sigma$-compatible Krein product on $K$. First we can suppose without loss of generality that $\rho$ is irreducible. Then we can use $\rho$ to tranfer $\times_\sigma$ from $\CC l(V)$ to an antilinear antiautomormphism $A\mapsto A^{\times_\sigma}$ of $End(K)$. Such an antiautomorphism is necessarily the adjoint operation for some non-degenerate hermitian form $(.,.)^K_\sigma$. To see this, use the simplicity of $End(K)$ to pick a $\beta$ such that $A^{\times_\sigma}=\beta A^\dagger \beta^{-1}$ where $A^\dagger$ is the adjoint for some scalar product $\bra .,.\ket$. Then it  is a simple matter to verify that $(.,.)_\sigma^K:=\bra .,\beta^{-1} .\ket$ has the required property.

Since this procedure is essentially in \cite{rob}, we will  call  it \emph{Robinson's transfert principle}.

\begin{propo}\label{rob1} There exists a non-degenerate hermitian form $(.,.)_\sigma^K$ on $K$ such that $\rho$ is $\sigma$-compatible, and when $\rho$ is irreducible it is unique up to multiplication by a non-zero real number.
\end{propo}

We haven't proven the uniqueness part yet. It is a simple consequence of Riesz'  representation theorem, but we state it here as a separate lemma  for future reference. The proof is left to the reader.

\begin{lemma}\label{unikrein} Let $(K,(.,.))$ be a space equipped with a non-degenerate hermitian form, and let $(.,.)'$ be a hermitian form such that for all $\psi,\phi\in K$, and for all $A\in B(K)$, one has $(A\psi,\phi)'=(\psi,A^\times \phi)'$, where $A^\times$ is the adjoint of $A$ for $(.,.)$. Then there exists $\lambda\in \RR$ such that $(.,.)'=\lambda (.,.)$.
\end{lemma}



 
Every minimal left ideal  of $\CC l(V)$ is an irreducible module for the representation of $\CC l(V)$ by left multiplication. Such a module is of the form $S_e:=\CC l(V)e$ where $e$ is a primitive idempotent of $\CC l(V)$, that is an idempotent which cannot be decomposed as a sum of two others. The module $S_e$ is sometimes called an \emph{algebraic spinor module}, and its elements are  {\it algebraic spinors}. It  is then  natural to restrict $(.,.)_\sigma$ to $S_e$ and ask if it defines a  $\sigma$-compatible Krein product. In fact $\sigma$-compatibility is immediate, but the restriction of $(.,.)_\sigma$  can be degenerate. Of course when $Q_\sigma$ is positive definite everything works fine : $(.,.)_\sigma$ is positive definite and its restriction to every algebraic spinor module is a $\sigma$-compatible scalar product. We will use this fact in the following proposition, which is implicit in \cite{rob}.

\begin{remark}
The positive definite case is  all we really need, but we can wonder what happens in general. It turns out that we just need to compute $ee^{\times_\sigma}$: either it is zero, in which case $(.,.)_\sigma$ vanishes on $S_e$, or it is not, in which case  $(.,.)_\sigma$ restricts to a 
$\sigma$-compatible Krein product on $S_e$. Moreover, is this latter case, there exists a unique primitive idempotent $f\in \CC l(V)$ such that $f^{\times_\sigma}=f$ and $S_e=S_f$. The proof of this claim is  given in appendix \ref{kreinalgspin} for the interested reader.
\end{remark}




\begin{propo}\label{propo7} (Robinson's alternative) Let $K$ be an irreducible spinor module and $(.,.)^K$ be a $\sigma$-compatible Krein spinor product.

If $Q_\sigma$ is positive definite, then $(.,.)^K$ is definite, if $Q_\sigma$ is not positive definite, then  $(.,.)^K$ is neutral.
\end{propo}
\begin{demo} Let us call $\rho : \CC l(V)\rightarrow End(K)$ the representation homomorphism.

If $Q_\sigma$ is not positive definite, let $e_i\in V_\sigma$ be such that $e_i^2=-1$, set $U=\rho(e_i)$ and use lemma \ref{neutralemma}.

If $Q_\sigma$ is positive definite, then let $U : K\rightarrow S_e$, where $S_e$ is a minimal left ideal in $\CC l(V)$,  be a Clifford-module isomorphism, which we know exists by irreducibility of $K$. By Garling's alternative, $(.,.)_\sigma$ is positive definite, hence its restriction to $S_e$ also is. We now transport $(.,.)_\sigma$ to $K$ thanks to $U$ by the formula $(\psi,\eta)_U:=(U\psi,U\eta)_\sigma$, for all $\psi,\eta\in K$. Let us check that $(.,.)_U$ is $\sigma$-compatible:

\bea
\forall a\in\CC l(V), (\rho(a)\psi,\eta)_U&=&(U\rho(a)\psi,U\eta)_\sigma\cr
&=&(L_aU\psi,U\eta)_\sigma\mbox{ by the intertwining property of }U\cr
&=&(U\psi,L_{a^{\times_\sigma}}U\eta)_\sigma\cr
&=&(U\psi,U\rho(a^{\times_\sigma})\eta)_\sigma\cr
&=&(\psi,\rho(a^{\times_\sigma})\eta)_U\nonumber
\eea

Now $(\psi,\psi)_U=0\Rightarrow (U\psi,U\psi)_\sigma=0\Rightarrow U\psi=0\Rightarrow \psi=0$, hence $(.,.)_U$ is definite, therefore $(.,.)^K$ is by lemma \ref{unikrein}.
\end{demo}

Of course multiplication by $-1$ turns $\sigma$-compatible Krein products into $\sigma$-compatible Krein products. Putting together Robinson's alternative and proposition \ref{rob1}, we can therefore characterize the Euclidean signature of $Q_\sigma$ in the following way.

\begin{cor} If the form $Q_\sigma$ is positive definite then  there exists a $\sigma$-compatible scalar product on  every irreducible spinor module. If there exists at least one $\sigma$-compatible scalar product on at least one irreducible spinor module then $Q_\sigma$ is positive definite.
\end{cor}

We now turn our attention to the case where $\sigma$ is admissible.

\begin{lemma}\label{lemme5}  Let $\rho$ be a  $c$-compatible representation of $\CC l(V)$ on a Krein space $(K,(.,.))$. Let $\sigma=Ad_b\circ c$ be an admissible real structure. Let $x\in \CC l(V)$ and let $(.,.)_x=(.,\rho(x).)$. Then $(.,.)_x$ is  a $\sigma$-compatible Krein product on $K$ iff $x=x^\times$ and $x$ is proportional to $b^{-1}$.
\end{lemma}
\begin{demo}
The two properties are easily seen to be sufficient. Let us prove that they are necessary. First we note that since $\CC l(V)$ is simple and $\rho\not=0$, then $\rho$ is injective.  For $(.,.)_x$ to be sesquilinear $x$ must satisfy $\rho(x)=\rho(x)^\times$ which is equivalent to $x=x^\times$ since $\rho$ is injective and $c$-compatible. Now we note that $a^{\times_\sigma}=\sigma(a^T)=\sigma\circ c\circ c (a^T)=Ad_b(a^\times)=ba^\times b^{-1}$. 

We must have $\rho(a^{\times_\sigma})=\rho(a)^{\times_\sigma}=\rho(x)^{-1}\rho(a)^\times \rho(x)$ and  we obtain that $xb$ is a scalar by the injectivity of $\rho$.
%
\end{demo}



Now if $b$ is real then $b^\times=\pm b$ according to proposition \ref{admreal}, and we can take $x=b^{-1}$ or $x=ib^{-1}$ in the above lemma. Let us store the result we obtain as a proposition.

\begin{propo}\label{gprod} Let $\rho$ be a  $c$-compatible representation of $\CC l(V)$ on a Krein space $(K,(.,.))$. Let $\sigma=Ad_b\circ c$ be an admissible real structure \emph{with $b$ real}.
\begin{enumerate}
\item If $b^\times=b$, $(.,.)_b:=(.,\rho(b)^{-1}.)$ is a $\sigma$-compatible Krein product on $K$.
\item If $b^\times=-b$, $(.,.)_b:=(.,i\rho(b)^{-1}.)$ is a $\sigma$-compatible Krein product on $K$.
\end{enumerate}
\end{propo}

\subsection{Characterization of the Lorentz and anti-Lorentz signatures}\label{sec27}
First we fix some terminology. If $Q$ has the Lorentz signature, the open light cone $C$ of $V$ is the non-convex   cone of those $v$ such that $Q(v)<0$. It consists of two connected components which we  arbitrarily call $C_+$ and $C_-$. If $Q$ has the anti-Lorentz signature, the open light cone $C$ is defined by $Q(v)>0$ and $C_+/C_-$ are defined accordingly.

Let $L$ be a subspace in $V$ and $s_L$ be the orthogonal symmetry with respect to $L$. That is, if $V=L\oplus W$ is an orthogonal decomposition, then $s_L$ is the identity on $L$ and minus the identity on $W$. 

Our starting point will be the following evident observation:

\begin{itemize}
\item $Q$ has the Lorentz signature iff there exists a line $L$ such that $B(-s_L(.),.)$ has Euclidean signature.
\item $Q$ has the anti-Lorentz signature iff there exists a line $L$ such that $B(s_L(.),.)$ has Euclidean signature.
\end{itemize}

Note that in both cases the line will belong to the open light cone. 

Now all we have to do is to translate this observation in terms of admissible real structures, and use the characterization of Euclidean signature we arrived at in the previous subsection. In order to do this we first review some basic facts about the implementation of orthogonal symmetries with the adjoint action of vectors in Clifford algebras. Special care is needed about the signs.

Let $v$ be a non-zero vector in $V$. Then $v$ is invertible and $v^{-1}={1\over Q(v)}v$. The adjoint action   of $v$ on a vector $u\in V$  satisfies $Ad_v(v)=vvv^{-1}=v$ and if $w\perp v$, $Ad_v(w)=vwv^{-1}=-wvv^{-1}=-w$. Thus $Ad_v=s_L$ where $L$ is the line $\RR v$.

Now notice that for any $u\in V$, $Ad_\omega(u)=\omega u\omega^{-1}=-u\omega\omega^{-1}=-u$ since $u$ is odd. Thus $Ad_{\omega v}=Ad_\omega \circ Ad_v=-s_L$.

Moreover, since $c(\omega)=\omega$ and $c(v)=v$, the real structures $\sigma_v:=Ad_v\circ c$ and $\sigma_{\omega v}:=Ad_{\omega v}\circ c$ are admissible. We have thus obtained the following:

\begin{lemma}\label{lem6} Let $L$ be a line in $V$. Then
\begin{itemize}
\item The unique admissible real structure which restricts to $s_L$ is $\sigma_v:=Ad_v\circ c$, where  $v\in L$ is a non-zero vector.
\item The unique admissible real structure which restricts to $-s_L$ is $\sigma_{\omega v}:=Ad_{\omega v}\circ c$, where  $v\in L$ is a non-zero vector.
\end{itemize}
\end{lemma}


Putting together what we have learned so far, we obtain:

\begin{theorem}\label{maintheorem} Let $V$ be a finite dimensional real vector space of even dimension with non-degenerate quadratic form $Q$. Let $\rho$ be a $c$-compatible irreducible representation on the Krein space $(K,(.,.))$. Then:
\begin{enumerate}
\item The signature of $Q$ is anti-Lorentzian iff there exists $v\in V$ such that $(.,.)_v:=(.,\rho(v)^{-1}.)$ is definite.  
\item If $n=2$ or $n=6$ modulo $8$, then the signature of $Q$ is  Lorentzian iff there exists $v\in V$ such that $(.,.)_{v}:=(.,\rho(\omega v)^{-1}.)$ is definite.  
\item If $n=0$ or $n=4$ modulo $8$, then the signature is  Lorentzian iff there exists $v\in V$ such that $(.,.)_{v}:=(.,i\rho(\omega v)^{-1}.)$ is definite.
\end{enumerate}

In every case,   $(.,.)_v$ is definite $\Leftrightarrow v\in C$, and there exists $\lambda=\pm 1$ such that $\lambda(.,.)_v$ positive definite $\Leftrightarrow v\in C_+$ and $\lambda(.,.)_v$ negative definite $\Leftrightarrow v\in C_-$.
\end{theorem}
\begin{demo}
\begin{enumerate}
\item We know that $Q$ has anti-Lorentz signature iff there exists a line $L$ such that $B(s_L(.),.)$ is positive definite. Observe that $B(s_L(.),.)=Q_\sigma$,  where $\sigma=Ad_v\circ c$, for any non-zero $v\in L$. We also know that $(.,.)_{v}$ is a $\sigma$-compatible Krein product on $K$ by proposition \ref{gprod}, and we conclude by Robinson's alternative.
\item Case 2 and 3 are similar to case 1 with the extra complication that $b=\omega v$ satisfies $b^\times=b$ for $n=2,6$ mod $8$ and $b^\times=-b$ for $n=0,4$ mod $8$.
\end{enumerate}
It is obvious that a line $L$ is such that $B(s_L(.),.)$ is positive definite iff $L\setminus\{0\}$ lies inside the light cone. Moreover multiplication by $-1$ exchanges $C_+$ with $C_-$ and $(.,.)_v$ with $-(.,.)_v$, so that it only remains to show that $(\psi,\psi)_v$ keeps a constant sign for a fixed non-zero $\psi$ and $v$ varying continuously in one component of $C$, which is immediate by continuity of $\rho$. 
\end{demo}

Cases 2 and 3 can be put together by saying that the signature is Lorentzian iff there  exists $v\in V$ such that $(.,.)_{v}:=(\chi\rho(v).,.)$ is definite, where $\chi$ is chirality operator of noncommutative geometry.

\begin{remark} Even if we could simplify the above theorem by taking $\rho(v)$ instead of $\rho(v)^{-1}$ and so on, since they differ only by a non-zero constant, it could prove misleading in the noncommutative situation. We will come back to this issue later.
\end{remark}
\smallbreak




We will not delve much into such matters, but the above characterization of the Lorentzian/antilorentzian signature can be generalized easily. If $Q$ is of signature $(p,q)$ and $p$ is odd, then we can find $p$ vectors $v_1,\ldots,v_p$ such that  $\sigma=Ad_g\circ c$, with $g=v_1\ldots v_q$, is an Euclidean real structure. The converse is also true, and these vectors will then automatically satisfy $v_i^2>0$. If $p$ is even we have to take $g=v_1\ldots v_q$ and we will have $v_i^2<0$. The problem now is to characterize the elements of the Clifford group which are of this form. For this we can use the canonical isomorphism $\Theta : \Lambda V\rightarrow \CC l(V)$. Thus we need to consider two different classes of quadratic forms: even and odd ones, meaning that $p$ and $q$ are both even or both odd, according to the case.  We can then say that:

\begin{itemize}
\item If $Q$ is odd, then $p$ is the only integer such that $\exists g\in \Theta(\Lambda^p V)$, $Ad_g\circ c$ is an Euclidean real structure.
\item If $Q$ is even then $q$  is the only integer such that $\exists g\in \Theta(\Lambda^q V)$, $Ad_g\circ c$ is an Euclidean real structure.
\end{itemize}

We could then translate these properties in terms of the Krein structures on spinor modules. However, in order to generalize to noncommutative geometry, we would have to use noncommutative $p$-forms or $q$-forms instead of elements of the exterior algebra over $V$, and this raises the issue of junk forms when $p$ (resp. $q$) is $>1$. We leave this question open for the moment, but we note that in the Lorentz case we have used the trick  of composing with the chirality in order to stick with 1-forms. It remains to see, using a more general approach, if this trick is legal. This is the main reason, besides simplicity, why we will mostly deal with the antilorentzian signature in this paper.

\subsection{Some Krein positive operators on spinor modules}\label{sec26}

In this subsection the quadratic form $Q$ is supposed to be antilorentzian. We consider  a $c$-compatible irreducible representation $\rho$ of $\CC l(V)$ into a Krein space $(K,(.,.))$, which we decompose as $K=K_L\oplus K_R$ where $K_L$ (resp. $K_R$) is the $-1$-eigenspace (resp. $+1$-eigenspace) of  the chirality operator $\chi$. It is immediately seen by computing $(\psi,\chi^2\phi)$ and using $\chi^\times=-\chi$ that $K_L$ and $K_R$ are self-orthogonal for the Krein product $(.,.)$.

An operator $A$ on $K$ will be said to be \emph{Krein-positive} iff $(\psi,A\psi)>0$ for all non-zero $\psi\in K$. We gather here some observations which will be useful later. For brevity, we do not write explicitly the representation $\rho$.

We already know that if $e$ is a timelike vector, $(.,.)_e$ is definite, and we will suppose here that the definition of the future light-cone and the Krein product are compatible in the sense that $e$ is future-directed iff it is Krein-positive. The next lemma indicates what happens if we take a spacelike vector instead.


\begin{lemma}\label{signspat} Let $w$ be a spacelike non-zero element of $V$ and let $(\psi,\psi)_w:=(\psi,w\psi)$. Then $K_L$ and $K_R$ are orhogonal to each other for $(.,.)_w$. If $n>2$, the form $(.,.)_w$ is neutral on each half-spinor module $K_L$ and $K_R$.
\end{lemma}
\begin{demo}
First, since multiplication by $w$ changes chirality, it is obvious that $K_L$ and $K_R$ are orthogonal to each other for $(.,.)_w$ since they are self-orthogonal for $(.,.)$. If $n>2$ we can find a spacelike vector $e$ such that $e^2=-1$ and $B(e,w)=0$. Let $U_e : K_L\rightarrow K_R$  be the map $\psi\mapsto e\psi$. Since $(e\psi,e\psi)_w=(e\psi,we\psi)=-(e\psi,ew\psi)=(\psi,\psi)_w$, $U_e$ is an isometry for $(.,.)_w$. Thus $(.,.)_w$, which is globally neutral, must have the same signature on $K_R$ and $K_L$, hence the result.
\end{demo}

If $n=2$ the form $(.,.)_w$  has opposite signature on the $1$-dimensional half-spinor modules.

There are other Krein-positive operators besides future-directed vectors.  An important example is given in the following lemma.

\begin{lemma}\label{lem8} Let $u,v\in V$. If $n>2$, $u+\chi v$ is Krein positive iff $u+v$ and $u-v$ are both timelike and future-directed.
\end{lemma}
\begin{demo}
First, using the fact that $K_L$ and $K_R$ are self-orthogonal we have
$$(\psi,(u+\chi v)\psi)=(\psi_L,(u- v)\psi_L)+(\psi_R,(u+ v)\psi_R)$$
with $\psi=\psi_L+\psi_R$, $\psi_{L/R}\in K_{L/R}$. If $u\pm v$ are timelike and future-directed then $(.,(u\pm v).)$ are both scalar products, hence the two summands are positive. Conversely, if $u+v$ is timelike and past-directed, we can take $\psi_L=0$ and $\psi_R\not=0$ to obtain a negative result. Similarly $u-v$ cannot be past-directed. If $u+v$ is spacelike we can use lemma \ref{signspat} to see that there exists a $\psi_R$ such that $(\psi_R,(u+v)\psi_R)<0$. A similar conclusion holds if $u-v$ is spacelike. 
\end{demo}

If $n>2$, a general $\times$-self-adjoint odd element of $End(K)$ will have the form $u+\chi v+r$, where $u,v\in V$ and the remainder $r$ is a linear combination of $e_I$, $I\subset\{1;\ldots;n\}$, $1<|I|<n-1$, with $(e_i)_{1\le i\le n}$ is some pseudo-orthogonal basis. Let us suppose that $u+\chi v+r$ is Krein-positive. Let us observe that if $a,b$ are Krein-positive then $\tau_n(ab)>0$. Indeed, for all non-zero $\psi\in K$, $(ab\psi,b\psi)>0$ by the Krein-positivity of $a$ and the invertibility of $b$. Hence $ab$ is a positive operator for the scalar product $(.,.)_b$, thus it must have a positive normalized trace.
Using this principle with $a=u+\chi v+r$ and $b$ a future-directed timelike vector, we see that $\tau_n(ub+\chi vb+rb)>0$. However $\tau_n(\chi vb)=\tau_n(rb)=0$ since $\chi v$ and $r$ are both orthogonal to $b$ with respect to the $c$-product. Hence we obtain that $\tau_n(ub)=B(u,b)>0$ for all future-directed timelike vectors $b$, which proves that $u$ is itself future-directed timelike. We state this result as a lemma.

\begin{lemma}\label{lem9} If $u+\chi v+r$ is Krein-positive, then $u$ is timelike and future-directed.
\end{lemma}

As an example of such a Krein-positive operator with a non-zero $r$, one can consider a pseudo-orthonormal basis $(e_1,\ldots,e_4)$ of $\RR^{1,3}$ and let $p=e_1+ite_2e_3$ for $t\in ]-1;1[$.

\subsection{Real structure and charge conjugation. KO-dimension tables}

In noncommutative geometry, the real structure is defined as an antilinear operator on the Krein space which is the abstract substitute of the space of sections of the spinor bundle. When it is viewed in this way, it is sometimes called the \emph{charge conjugation}, and we will adopt this terminology in order to distinguish it from its Clifford counterpart. To begin this section let us see how the two notions relate locally.

\begin{propo}\label{generalWR} Let $\rho : \CC l(V)\rightarrow End(K)$ be an irreducible $c$-compatible representation. Then there exists an antilinear operator ${\mathcal C} : K\rightarrow K$ implementing $c$, i.e. such that
$$\rho(c(a))={\mathcal C}\rho(a){\mathcal C}^{-1}$$
 for all $a\in \CC l(V)$. This operator can be chosen to satisfy ${\mathcal C}^2={\tilde \epsilon}$ and ${\mathcal C}^\times {\mathcal C}={\tilde \kappa}$, with ${\tilde \epsilon}$ and ${\tilde \kappa}$ some signs, and in that case  is  unique up to multiplication by $e^{i\theta}$, $\theta\in\RR$.

Let $\sigma=Ad_b\circ c$ be an  admissible real structure, with $b$ in the Clifford group such that $b^2=\lambda=\pm 1$, $b^\times=\lambda' b$, $\lambda'=\pm 1$, and $c(b)=b$. Let $B=\rho(b)$. Then $\sigma$ is implemented on $K$ by  ${\mathcal C}_\sigma=B{\mathcal C}={\mathcal C}B$.

We have ${\mathcal C}_\sigma^2=\lambda {\tilde \epsilon}$, ${\mathcal C}_\sigma^\times {\mathcal C}_\sigma=\lambda\lambda'{\tilde \kappa}$, and if $( .,.)_\sigma$ is
 a $\sigma$-compatible Krein product on $K$, then $({\mathcal C}_\sigma)^{\times_\sigma}=\lambda{\mathcal C}_\sigma^\times$ and $({\cal C}_\sigma)^{\times_\sigma}{\cal C}_\sigma=\lambda'{\cal C}^\times {\cal C}$.

\end{propo}
\begin{demo}
Fix any basis of $K$ and denote by $c.c. : \psi\mapsto\bar\psi$ the complex conjugation of coordinates with respect to this basis. We also denote by $c.c. : A\mapsto \bar A$ the antilinear involution on operators induced by complex conjugation of coordinates, that is : $\bar A\psi=\overline{A\bar \psi}$.  Let us write  $\tilde c=\rho \circ c\circ \rho^{-1}$. Since $\tilde c$ is an antilinear involution of $End(K)$, $\tilde c \circ c.c.$ is a linear (involutive) automorphism and is thus of the form $Ad_C$ for some $C\in End(K)$. We let ${\cal C}=C\circ c.c.$ and we obtain  that $\rho(c(a))={\mathcal C}\rho(a){\mathcal C}^{-1}$.


The fact that $c$ is an involution translates as ${\mathcal C}^2={\tilde \epsilon}$, with ${\tilde \epsilon}$ a constant. Let $\psi\in K$ be such that ${\mathcal C}\psi\not=0$. Then ${\mathcal C}^3\psi={\tilde \epsilon} {\mathcal C}\psi={\mathcal C}({\tilde \epsilon}\psi)$, which shows that ${\tilde \epsilon}$ is real. From the fact that $\rho$ is $c$-compatible we obtain that $\rho(c(a))^\times=\rho(c(a)^\times)=\rho(c(a^\times))$ and it boils down to ${\mathcal C}{\mathcal C}^\times={\tilde \kappa}$, with ${\tilde \kappa}$ a constant, which must be real. Calculating ${\mathcal C}^2({\mathcal C}^\times)^2$ we find that ${\tilde \kappa}^2={\tilde \epsilon}^2$. Dividing ${\mathcal C}$ by $\sqrt{|{\tilde \kappa}|}$  we can suppose that ${\tilde \kappa}=\pm 1$ and ${\tilde \epsilon}=\pm 1$. It is now easy to prove the uniqueness up to a phase.  

The second part is an obvious consequence of proposition \ref{admreal}. We have ${\mathcal C}_\sigma^\times {\mathcal C}_\sigma={\mathcal C}^\times B^\times B{\mathcal C}=\lambda\lambda'{\mathcal C}^\times {\mathcal C}=\lambda\lambda'{\tilde \kappa}$. If $\lambda=1$ then $({\mathcal C}_\sigma)^{\times_\sigma}=B^{-1}{\mathcal C}_\sigma^\times B={\mathcal C}_\sigma^\times$ (see proposition \ref{gprod}), and if $\lambda=-1$, then $({\mathcal C}_\sigma)^{\times_\sigma}=(iB)^{-1}{\mathcal C}_\sigma^\times (iB)=-{\mathcal C}_\sigma^\times$.

Hence we have $({\cal C}_\sigma)^{\times_\sigma}{\cal C}_\sigma=\lambda{\cal C}_\sigma^\times{\cal C}_\sigma=\lambda {\cal C}^\times B^\times B{\cal C}=\lambda\lambda'{\cal C}^\times B^2{\cal C}=\lambda'{\cal C}^\times {\cal C}$.
\end{demo}

If the representation is not irreducible, we can break it up into irreducible parts and use the proposition on each one of them. The signs $\tilde \epsilon$ and $\tilde \kappa$ can still be defined globally since they do not depend on the representation, but the phase $\theta$ can vary from block to block.

Let us use   proposition \ref{generalWR} to see how the KO table of signs change when we perform a Wick rotation. 

If we start from  an antilorentzian quadratic form $Q$, a Wick rotation to Euclidean signature corresponds to $b=v$ with $v^2=1$, $v^\times=v$, hence $\lambda'=\lambda=1$.

In the Lorentzian case $b=\omega v$ with $v^2=-1$, hence $b^2=-\omega^2 v^2=(-1)^{{n\over 2}+1}$, and $b^\times=v^\times\omega^\times=(-1)^{n\over 2}v\omega=(-1)^{{n\over 2}+1}b$. Hence one  has $\lambda=\lambda'=(-1)^{{n\over 2}+1}$ in this case. 

We introduce the hopefully obvious notations ${\cal C}_L$, $\chi_L$, ${\cal C}_E$, $\chi_E$, ${\cal C}_{AL}$, $\chi_{AL}$. Then ${\cal C}_E={\cal C}_{AL} \rho(v)={\cal C}_{L}\rho(\omega v) $. According to the above proposition we then have:

$$ {\cal C}_E^2={\cal C}_{AL}^2 ;\   {\cal C}_E^*{\cal C}_E={\cal C}_{AL}^\times{\cal C}_{AL} $$

when Wick rotating from   antilorentzian to   Euclidean signature, but

$${\cal C}_E^2=(-1)^{{n\over 2}+1}{\cal C}_L^2;\  {\cal C}_E^*{\cal C}_E=(-1)^{{n\over 2}+1} {\cal C}_{L}^\times{\cal C}_{L} $$

when Wick rotating from   Lorentzian to   Euclidean signature. Here we have used the notation $*$ instead of $\times_\sigma$ since we are dealing with a Hilbert adjoint.

When going from antilorentzian/Lorentzian to Euclidean signature, the volume element changes in the following way : $\omega_{AL/L}=e_1\ldots e_{n}\mapsto \omega_E:=u_\sigma(\omega_{AL/L})=i^q\omega_{AL/L}$. But the Euclidean chirality element is $\chi_E=(-i)^{n\over 2}\rho(\omega_E)$ whereas $\chi_{AL/L}=(-i)^{{n\over 2}+q}\rho(\omega_{AL/L})$. Hence the two expressions only differ by a minus sign:

$$\chi_E=-\chi_L, \chi_E=-\chi_{AL}$$

Finally we note that in both cases ${\cal C}_E$ has an additional $-1$ factor in its commutation relation with $\chi_E$ with respect to the one of ${\cal C}_{AL/L}$. 

Summarizing we have:

\begin{itemize}
\item Antilorentzian $\rightarrow$ Euclidean : $\tilde\epsilon\mapsto \tilde\epsilon$, $\tilde\kappa\mapsto\tilde\kappa$, $\tilde\epsilon''\mapsto -\tilde\epsilon''$
\item Lorentzian $\rightarrow$ Euclidean : $\tilde\epsilon\mapsto (-1)^{{n\over 2}+1}\tilde\epsilon$, $\tilde\kappa\mapsto(-1)^{{n\over 2}+1}\tilde\kappa$, $\tilde\epsilon''\mapsto -\tilde\epsilon''$
\end{itemize}

These rules permit to fill the   Lorentzian and antilorentzian KO table of signs from the Euclidean one, using the fact that the metric dimension modulo 8 is preserved by a Wick rotation (the KO dimension, of course, is not). 

The real Clifford algebra $Cl(V,Q)$ can be directly recovered as the algebra of all $a\in \CC l(V)$ such that  $[{\mathcal C},\rho(a)]=0$.  In noncommutative geometry, one generally  uses the operator $J$ which selects $Cl(V,-Q)$ (as originally in \cite{ncgreal}). To pass from one to the other is easy:

$$J=\chi{\mathcal C}$$

We call ${\cal C}$ the \emph{ungraded charge conjugation} and $J$ the \emph{graded charge conjugation} operator. In this work it is generally more convenient to use ${\cal C}$, but we stick to the traditional $J$ in definitions and statements of theorem.

The signs $\epsilon,\epsilon',\epsilon''$ which are given in the definition of spectral triples correspond to the commutation rules of $J$ with $D$ and $\chi$ and to the square of $J$. In the Lorentzian/antilorentzian context there is also a sign arising from $J^\times J$.

When we pass from the ${\mathcal C}$ convention (signs with tildes) to the $J$ convention (signs without tildes)  and vice versa we must use the following rules:

\begin{itemize}
\item $\tilde\epsilon''=\epsilon''$ ($J$ and ${\cal C}$ have the same commutation rules with $\chi$).
\item $\tilde\epsilon'=-\epsilon'$ ($J$ and $D$ have a commutation sign opposite to the one of ${\cal C}$ and $D$). In even dimension $\epsilon'=1$ and $\tilde\epsilon'=-1$ always. Hence we do not let these signs appear in the table.
\item $\tilde\epsilon=\epsilon''\epsilon$ (because $J^2=\chi{\mathcal C}\chi{\mathcal C}=\epsilon''{\mathcal C}^2$).
\item In the Lorentzian/antilorentzian case : $\tilde\kappa=-\kappa$ (since $J^\times J=(\chi {\mathcal C})^\times \chi{\mathcal C}=-{\mathcal C}^\times{\mathcal C}$), in the Euclidean case, $\tilde\kappa=\kappa$.
\end{itemize}

\begin{table}
\begin{center}
\begin{tabular}{ccccc}
KO dim $=n=p-q\ [8]$ & 0 & 2 & 4 & 6 \\
$J^2=\epsilon$ & 1 & -1 & -1 & 1 \\
$J\chi=\epsilon''\chi J$, ${\cal C}\chi=\epsilon''\chi {\cal C}$ & 1 & -1 & 1 & -1 \\${\cal C}^2=\tilde\epsilon$ & 1 & 1 & -1 & -1 \\
$J^\times J=\kappa$ & 1 & 1 & 1 & 1 \\
${\cal C}^\times {\cal C}=\tilde\kappa$ & 1 & 1 & 1 & 1 
\end{tabular}
\caption[smallcaption]{Euclidean table of  KO signs}
\label{KOeuc}
\end{center}
\end{table}

\begin{table}
\begin{center}
\begin{tabular}{ccccc}
$n \ [8]$ & 0 & 2 & 4 & 6 \\
KO dim $=p-q\ [8]$ & 2 & 0 & 6 & 4 \\
$\epsilon$ & -1 & 1 & 1 & -1 \\
$\epsilon''=\tilde\epsilon''$ & -1 & 1 & -1 & 1 \\
$\tilde\epsilon=\epsilon''\epsilon$ & 1 & 1 & -1 & -1 \\
$\kappa$ & -1 & -1 & -1 & -1 \\
$\tilde\kappa$ & 1 & 1 & 1 & 1 
\end{tabular}
\caption[smallcaption]{Antilorentzian table of  KO signs.}
\label{KOantilor}
\end{center}
\end{table}

\begin{table}
\begin{center}
\begin{tabular}{ccccc}
$n$ $[8]$& 0 & 2 & 4 & 6 \\
KO dim $=p-q\ [8]$ & 6 & 0 & 2 & 4 \\
$\epsilon$ & 1 & 1 & -1 & -1 \\
$\epsilon''=\tilde\epsilon''$ & -1 & 1 & -1 & 1 \\
$\tilde\epsilon=\epsilon''\epsilon$ & -1 & 1 & 1 & -1 \\
$\kappa$ & 1 & -1 & 1 & -1 \\
$\tilde\kappa$ & -1 & 1 & -1 & 1 
\end{tabular}
\caption[smallcaption]{Lorentzian  table of KO signs.}
\label{KOlor}
\end{center}
\end{table}


\subsection{Examples}

In this section, which is optional, we derive the compatible Krein products on spinors in some usual representations. We illustrate how the light cone appears in each of the three different cases of theorem \ref{maintheorem} in concrete computations with gamma matrices. Let us fix the notations. In each case we take $V=\RR^n$, with a pseudo-orthonormal basis $(e^\mu)_{\mu=0,\ldots,n-1}$. We have raised the indices to be consistent with physicists conventions. We will have $B(e^0,e^0)=-1$ in the Lorentz case and $B(e^0,e^0)=1$ in the anti-Lorentz case, and $e^0$ will belong to $C^+$ in both cases.

We will use a representation (for instance Dirac or Majorana\footnote{In fact the Dirac representation is more natural, since the Dirac basis is at the same time orthonormal for the scalar product $\bra .,.\ket_{e^0}$ and pseudo-orthonormal for the Krein-product.}) in which the gamma matrices $\gamma^\mu:=\rho(e^\mu)$  satisfy $(\gamma^\mu)^\dagger=\pm\gamma^\mu$ where the sign is the same as in $(\gamma^\mu)^2=\pm 1$, and $\dagger$ is the adjoint relatively to the canonical scalar product on $K=\CC^N$, $N=2^{n/2}$.

We will consider a $c$-compatible Krein product on $K$ denoted by $(.,.)$. Its adjunction operation is denoted $A\mapsto A^\times$. We denote by $\beta$ the matrix of the form $(.,.)$ in the canonical basis, that is, $(.,.)=\bra .,\beta .\ket$ where $\bra .,.\ket$ is the canonical scalar product. We know that $(.,.)$ is determined up to a constant, which we reduce to a sign  by requiring that $\beta^2=1$.

The $c$-compatibility of $\rho$ is equivalent to the requirement that $(\gamma^\mu)^\times=\gamma^\mu$ for all $\mu$, which is in turn equivalent to $\beta(\gamma^\mu)^\dagger\beta^{-1}=\gamma^\mu$. We see then that :

\begin{itemize}
\item In the anti-Lorentz case, $\beta$ is   determined up to a sign by :
\be
\beta\gamma^0\beta^{-1}=\gamma^0,\quad \beta\gamma^k\beta^{-1}=-\gamma^k, k=1,\ldots,n-1,\quad \beta^\dagger=\beta, \beta^2=1\label{eq1}
\ee
\item In the Lorentz case,  $\beta$ is   determined up to a sign by :
\be
\beta\gamma^0\beta^{-1}=-\gamma^0,\quad \beta\gamma^k\beta^{-1}=\gamma^k, k=1,\ldots,n-1,\quad \beta^\dagger=\beta,\beta^2=1\label{eq2}
\ee
\end{itemize}

The obvious solution to (\ref{eq1}) is $\beta=\pm\gamma^0$. The more or less obvious solution to (\ref{eq2}) is $\beta=\pm\gamma^0\rho(\omega)$ if $n=2,6$ $[8]$ and $\beta=\pm i\gamma^0\rho(\omega)$ if $n=0,4$ $[8]$. In other words  in either case $\beta=\pm\gamma^0\chi$ where $\chi$ is the chirality operator (the ``$\gamma^5$ matrix'' in this context).

Now the part about the light cone in theorem \ref{maintheorem} tells us in the anti-Lorentz case that $v\in C$ if and only if $\beta\rho(v)$ is a definite hermitian matrix. In the Lorentz case  it is $\beta\chi\rho(v)$ which must be definite, but since the two $\chi$ cancel, we obtain   in both cases that $v\in C\Leftrightarrow \gamma^0\rho(v)$ is a definite hermitian matrix.

\section{Global constructions}

Let $(M,g)$ be a connected orientable semi-Riemannian manifold of dimension $n$ and signature $(p,q)$. In this section the notion of spacelike (resp. timelike) will be associated with the positive (resp. negative) index of inertia. Note that this clashes with the relativistic convention in the antilorentzian case ! We see no way of avoiding this problem. We will warn the reader when a confusion could arise. We will revert to the usual convention in the following sections when we are exclusively concerned with the antilorentzian/Lorentzian cases.

Given the Clifford bundle $Cl(M,g)$  and its complexification $\CC l(M)$ we now seek to make the constructions of the preceding section global. We start with real structures.

\subsection{Global real structures}

A  \emph{(global) real structure} $\sigma$ on $\CC l(M)$ is an involutive bundle map which is antilinear and respect products over each fibre. Since $\CC l(M)$ is given as the complexification of $Cl(M,g)$, we start with a given real structure $c$. We will be interested in real structures $\sigma$ commuting with $c$ and such that the metric $g(\sigma(.),.)$ is Euclidean, i.e. \emph{Euclidean real structures}. In view of proposition \ref{propo6} such a real structure $\sigma$ defines an orthogonal splitting  $TM=E_s\oplus E_t$ of the tangent bundle of $M$ into spacelike and timelike  subbundles, and the restriction of $\sigma$ to $TM$ is the orthogonal symmetry  with respect to $E_s$. Conversely, given such a splitting (which always exists,   see \cite{finster}, section 3), we can define the spatial orthogonal symmetry  globally, and extend it   by the universal property of the Clifford bundle. If we compose the result with $c$ we obtain an Euclidean real structure. Hence there is no obstruction to the existence of such objects.

However, $\sigma$ is given at each point $x\in M$ by $\sigma_x=Ad_{b_x}\circ c_x$ where $b_x$ in an element of the Clifford group at $x$. Hence globally we will have $\sigma=Ad_b\circ c$ where $b$ is a section of the Clifford group such that $Ad_b$ is smooth. But it is clear that there will generally be an obstruction   to the existence of a smooth $b$. To get a feel of what this obstruction might be let us consider the antilorentzian case. We know that  $\sigma$ defines an orthogonal splitting  $TM=E_s\oplus E_t$ of the tangent bundle, and in view of  lemma \ref{lem6},  $b$ itself is a section of $E_s$. Asking this section to be smooth   is asking for  the existence of a \emph{time}-orientation\footnote{Reader beware ! What we have called spacelike in this semi-Riemannian context is what is called timelike in Relativity when the mostly plus convention is adopted.}, that is a nonvanishing timelike (in the sense of relativity) vector field (\cite{beem}, p 5).

To deal with the general case we need a formal definition of space and time orientations in the   semi-Riemannian case\footnote{Surprisingly we could not locate any  in the literature.}. Recall that the kernel of a $p$-form $\omega_x$ is the subspace of $T_xM$ consisting of those vectors $v$ such that $\omega(v,.,\ldots,.)$ is the null $p-1$-form.

\begin{lemma} The following claims are equivalent.
\begin{enumerate}
\item\label{p1} There exists  a $p$-form $\omega$ such that $\ker(\omega_x)$ is a timelike $q$-dimensional subspace of $T_xM$ for all $x\in M$.
\item\label{p2} There exists a $p$-form $\omega$ such that for every linearly independent family $(v_1,\ldots,v_p)$ of spacelike tangent vectors at $x$, $\omega_x(v_1,\ldots,v_p)\not=0$.
\item\label{p3} Given any decomposition $TM=E_s\oplus E_t$ of the tangent bundle into the sum of a spacelike and timelike subbundle, there exists a non-vanishing top form $\omega^s$ of $E_s$.
\end{enumerate}
\end{lemma}
\begin{demo}
To see that (\ref{p1}) entails (\ref{p2}) consider any   family ${\cal F}=(v_1,\ldots,v_p)$ of spacelike tangent vectors at $x$ and let $S_x\subset T_xM$ be the linear subspace they span. Since $\ker\omega_x$ is timelike and $q$-dimensional, one has $\ker \omega_x\oplus S_x=T_xM$. Call $\pi_x$ the projection onto $S_x$ defined by this decomposition. It is immediate that $\omega_x(u_1,\ldots,u_p)=\lambda\det_{{\cal F}}(\pi_x(u_1),\ldots,\pi_x(u_p))$ for any vectors $u_1,\ldots,u_p$, where $\det_{{\cal F}}$ is the determinant in the basis ${\cal F}$ of $S_x$ and $\lambda\in\RR$. Clearly $\lambda\not=0$ since $\omega_x$ does not vanish identically. Hence $\omega_x(v_1,\ldots,v_p)\not=0$.

The proof that (\ref{p2}) entails (\ref{p3}) is immediate. To obtain (\ref{p1}) from (\ref {p3}), consider the projection $\pi$ on $E_s$ defined by $TM=E_s\oplus E_t$, and extend $\omega^s$ by the formula 
$$\omega_x(u_1,\ldots,u_p):=\omega^s_x(\pi(u_1),\ldots,\pi(u_p))$$
Then $\omega$ is a $p$-form on $TM$ and its kernel at $x$ is clearly $(E_t)_x$.
\end{demo}

A similar lemma clearly holds with timelike/spacelike reversed. When the properties stated in the lemma hold, we will say $(M,g)$ is \emph{space (respectively) time orientable}. A form on $TM$ with the two first properties will be called a space (resp. time) orientation.

If  $\omega'$ is a space orientation $p$-form on $M$, a linearly independent family $(v_1,\ldots,v_p)$ of spacelike vectors at $x$ will be said to be positively oriented iff $\omega_x'(v_1,\ldots,v_p)>0$. A similar definition can be given with time-orientations. Note that  space/time orientation forms actually provide more information than   the notion of orientation on families of spacelike/timelike vectors just defined. For example, suppose a Lorentz manifold is time-orientable and consider a $1$-form $\omega''$ with spacelike kernel. Then a timelike vector $v$ at $x$ will be said to be positively oriented, or \emph{future-directed} iff $\omega_x''(v)>0$. Using the musical isomorphism $\sharp$ provided by $g$ one can define a vector field $\xi=(\omega'')^\sharp$, which is timelike, since it is orthogonal to the distribution of spacelike subspaces $\ker(\omega'')$, and non-vanishing. This is indeed the usual definition of time-orientation for Lorentz manifolds, even if the notion of future/past directed timelike vectors would just require  the continuous choice of a half light-cone at each point of the manifold.

In the same vein  we say that $M$ is (totally) orientable if there exists a $n$-form $\omega$ such that $\omega_x$ has zero kernel (that is, $\omega$ is non-vanishing) for all $x$. 





Now the canonical isomorphism $\Theta_x : \Lambda T_x M\simeq Cl(T_x M,g_x)$ gives rises to a bundle isomorphism (see for instance \cite{LM}, prop. 3.5, chap. 2). We can compose it with the musical isomorphism and obtain a canonical isomorphism of vector bundles $\Theta : \Lambda T^* M\simeq Cl(M,g)$. If $\omega$ is a space orientation form, then at each $x$ we can decompose $T_xM$ into an orthogonal sum $\ker(\omega_x)\oplus S_x$, where $S_x$ is spacelike. If ${\cal B}=(v_1,\ldots,v_p)$ is an orthonormal basis of $S_x$, we obtain immediately that $\omega_x(u_1,\ldots,u_p)=\lambda_x \det_{\cal B}(\pi(u_1),\ldots,\pi(u_p))$ with $\pi$ the orthogonal projection on $S_x$. Hence $b_x:=\Theta_x(\omega_x)=\lambda_x v_1\ldots v_p$ is  seen to be an element of the Clifford group at $x$. Thus $b=\theta(\omega)$ is a smooth section of the Clifford group bundle. Conversely if $b$ is smooth section of the Clifford group bundle which is locally of the form $v_1\ldots v_p$, with $v_i$ spacelike, we easily see that $\theta^{-1}(b)$ is a spatial orientation. Similar consideration apply to time orientations.

We can use these observations to globalize the discussion at the end of subsection \ref{sec27}. Moreover since we have assumed $M$ to be   orientable, space orientations can be converted into time orientations,  saving us the trouble distinguishing the even and odd cases, and of translating the words timelike/spacelike when we revert to the traditional convention. We can thus answer our question about the obstruction to the definition of real structure through smooth sections of the Clifford group bundle in the following way:

\begin{itemize}
\item An orientable semi-Riemannian manifold admits Euclidean real structures defined through smooth sections of the Clifford group bundle iff it is time and space orientable.
\end{itemize}

Now recall that a time-oriented Lorentzian/antilorentzian manifold is called a   \emph{spacetime}. Hence an oriented Lorentzian/antilorentzian manifold admits Euclidean real structures defined through smooth sections of the Clifford bundle iff it is a \emph{spacetime}.


\subsection{Hermitian forms on the spinor bundle}

In this subsection we consider  a spin-c manifold $M$ with metric $g$ of signature $(p,q)$ and    spinor bundle $S$. We call $\rho$ the representation $\rho : \CC l(M)\rightarrow End(S)$.  We will often write $a.\Psi$ instead of $\rho(a)\Psi$ in order to simplify notations. Let $H : x\mapsto H_x$ be a smooth field of nowhere degenerate hermitian forms on $S$. We call this object a \emph{spinor metric}. For any $A\in End(S)$ we denote by $A^\times$ the map $x\mapsto A(x)^\times$, where $A(x)^\times$ is the adjoint relatively to $H(x)$. As in the local case, we say that $H$ is $c$-compatible if $\rho(a)^\times=\rho(a^\times)$ for all sections $a$ of $\CC l(M)$. This is equivalent to ask that $\rho(\xi)$ be self-adjoint for all vector fields $\xi\in \Gamma(TM)$. What we want now is to globalize Robinson's transfert principle.  The following theorem extends results in \cite{baum}.

\begin{theorem}\label{theorem2} There exists  a $c$-compatible spinor metric iff $(M,g)$ is time orientable when $p,q$ are even, and iff $(M,g)$ is space orientable when $p,q$ are odd.   
\end{theorem}
\begin{demo} We deal with the odd case only, the even case being completely similar.

Let us suppose that there exists a $c$-compatible spinor metric $H$. Consider a covering of $M$ by open sets $(U_\alpha)_{\alpha\in A}$ such that $TM$ and $S$ are trivial over $U_\alpha$, and a subordinate partition of unity $(f_\alpha)_{\alpha\in A}$. For each $x\in M$ we let $I_x$ be the finite set of indices $\alpha$ such that $f_\alpha(x)\not=0$. On each $U_\alpha$ let us choose a section $\psi_\alpha$ of the spinor bundle which is constant in some trivialization $S_{|U_\alpha}\simeq U_\alpha\times \CC^{2^{n/2}}$ and non-vanishing.

Let us define the $p$-form

\be
\omega(X_1,\ldots,X_p)=\sum_\alpha f_\alpha(x)H_x(\psi_\alpha,i^r\Theta(X_1\wedge\ldots\wedge X_s)\psi_\alpha)\nonumber
\ee

The integer $r$ is $[{p\over 2}]$. It ensures that $\omega$ is a real $p$-form. Indeed:

\bea
\Theta(X_1\wedge\ldots\wedge X_p)^\times&=&{1\over p!}\sum_{\sigma}\epsilon(\sigma)X_{\sigma(p)}\ldots X_{\sigma(1)}\cr
&=&{(-1)^{r}\over p!}\sum_{\tau}\epsilon(\tau)X_{\tau(1)}\ldots X_{\tau(p)}\cr
&=&(-1)^{r}\Theta(X_1\wedge\ldots\wedge X_p)\nonumber
\eea
where $(-1)^{r}$ is the signature of the reversal permutation $(1,\ldots,p)\mapsto (p,\ldots,1)$. Let $(e_1,\ldots,e_p,e_{p+1},\ldots,e_n)$ be a pseudo-orthonormal basis of $T_x M$ such that $e_1,\ldots,e_p$ are spacelike. Then we know  that $H_x(.,i^r(e_1\ldots e_p)^{-1}.)$ is definite. Since $(e_1\ldots e_p)^{-1}=e_p\ldots e_1=(-1)^re_1\ldots e_p$,  we have either:

$$\forall \alpha\in I_x, H_x(\psi_\alpha,i^r e_1\ldots e_p\psi_\alpha)>0$$

or

$$\forall \alpha\in I_x, H_x(\psi_\alpha,i^r e_1\ldots e_p \psi_\alpha)<0$$

In both cases we have $\omega_x(e_1,\ldots,e_p)\not=0$. Now, since $\omega$ is multilinear alternate, we have for any vectors   $u_1\ldots,u_p\in {\rm Span}(e_1,\ldots,e_p)$: 

\be
\omega_x(u_1,\ldots,u_p)=\det(u_i^j)_{1\le i,j\le p}\omega_x(e_1,\ldots,e_p)\label{eq10}
\ee
where $u_i^j$ is the $j$-th component of $u_i$ in the basis $(e_1,\ldots,e_p)$. Since $(e_1,\ldots,e_p)$ is any orthonormal family of spatial vectors, this shows that $\omega_x(u_1,\ldots,u_p)\not=0$ for any linearly independent family of spatial vectors. Thus $\omega$ is a space orientation.
 
Conversely, since $M$ is space-orientable, there exists a Euclidean real structure $\sigma=Ad_b\circ c$ where $b$ is a smooth section of the Clifford group bundle. On the Riemannian manifold $(M,g_\sigma)$ it is well-known that there exists a $\sigma$-compatible spinor metric. Let us call $\tilde H$ such a spinor metric. Then $H(.,.)=\tilde H(.,i^r b.)$ defines a $c$-compatible spinor metric for $M$.
\end{demo}

Of course a $c$-compatible spinor metric is unique only up to multiplication by a nonvanishing real function. 


In view of the above theorem we come to the conclusion that the  framework we will have a chance to   generalize to the noncommutative world  is that of space and time-orientable spin manifolds. We think   it   noteworthy that these purely mathematical considerations have led us to manifolds which are ``physics friendly'' in the sense that 1) we can put matter fields on them, 2) we can define an arrow of time, 3) particles will not change their chirality if transported around the universe. 

\subsection{The Dirac operator and $c$-compatibility}\label{sec33}

In the case of a manifold we have seen that we can express neatly the $c$-compatibility of a spinor metric by saying that vector fields are self-adjoint. Unfortunately this cannot be directly generalized in the noncommutative setting. It is therefore   important to have an alternative formulation which is better suited to this generalization. This formulation is quite simple: the Dirac operator has to be self-adjoint. We could close this section there, but we prefer to develop a little bit in order to be completely precise and also because it will be useful later on.

We will need to consider a spin manifold. Our definition of a spin structure on a spin-c manifold $(M,g,S)$ is the algebraic one, which is directly applicable in noncommutative geometry: it is antilinear bundle map $J : S\rightarrow S$ which satisfies $J\rho(a)J^{-1}=-\rho(c(a))$ and $J^2=\pm \id_S$. If $S$ is equipped with a $c$-compatible spinor metric $H$, we will also have $JJ^\times=\pm \id_S$ (see proposition \ref{generalWR}). We refer to \cite{schroder} for the equivalence of this definition and the traditional one involving the possibility of lifting the frame bundle to a spin bundle. 

By progressively enriching the structure we can define various kinds of connections on $S$. We will use the notation

\be
(\Psi,\Psi')_H=\int_MH_x(\Psi(x),\Psi'(x))\sqrt{|\det(g)|}dx\label{kreinprod}
\ee

We call this the \emph{Krein product} on spinor fields, a name which will be justified later on.

\begin{definition}\label{connections} Let $(M,g,S)$ be a spin-c manifold and $H$ be a spinor metric on $S$. Let $\nabla$ be a connection on $S$ .
\begin{enumerate}
\item If  $\nabla_X(a\cdot\Psi)=(\nabla_X^{LC}a)\cdot \Psi+a\cdot \nabla_X\Psi$,  for all sections $a$ of the Clifford bundle, vector fields $X$ and spinor fields $\Psi,\Psi'$ with compact support, then $\nabla$ is said to be a \emph{Clifford connection}. 
\item If $X\cdot (\Psi,\Psi')_H=(\nabla_X\Psi,\Psi')_H+(\Psi,\nabla_X\Psi')_H$, for all $X, \Psi,\Psi'$ as above, $\nabla$ is said to be \emph{metric}.
\item  If in addition $M$ is a spin manifold and  $J$ is a spin structure on it, then $\nabla$ is said to  \emph{preserves spin} if  
 $\nabla_X J=J\nabla_X$ for all vector field $X$.
\item If   $M$ is a spin manifold and is space and time oriented, $J$ is a spin structure on it, and $H$ is a $c$-compatible metric on $S$, then $\nabla$ is said to be \emph{a spin connection} if it is a Clifford connection which is metric and preserves spin.
\end{enumerate}
\end{definition}

 Let $x,y\in M$, let $\lambda$ be a curve joining these points. Let  $h_\lambda : S_x \rightarrow S_y$ be the parallel transport operator of $\nabla$ along  $\lambda$. Similarly, we denote by $h_\lambda^{LC}$ the parallel transport of the Levi-Civita connection along $\lambda$. Since this is an isometry from  $(T_xM, g_x)$ onto $(T_y M, g_y)$ we can consider its canonical extension $\tilde h_\lambda^{LC}$ which is an isomorphism between $\CC l(T_x M)$ and $\CC l(T_y M)$. The above infinitesimal definitions can be given   integrated forms.   More precisely, with the same hypotheses as in definition \ref{connections}, we can say that

\begin{enumerate}
\item $\nabla$ is a Clifford connection  iff for all $x,y,\lambda$, $h_\lambda$ intertwines the action of the Clifford algebras at $x$ and at $y$ on the spinors. More precisely, this means that  for all $a\in \CC l(T_x M)$ the following diagram commutes

$$\xymatrix{S_x\ar[r]^a\ar[d]_{h_\lambda}& S_x\ar[d]^{h_\lambda}\cr S_y\ar[r]^{\tilde h_\lambda^{LC}(a)}& S_y}$$
 
\item $\nabla$ is metric iff the parallel transport $h_\lambda$ is an isometry from $(S_x,H_x)$ onto $(S_y,H_y)$.
\item $\nabla$ preserves spin iff the following diagram commutes for all $x,y,\lambda$:

$$\xymatrix{S_x\ar[r]^{J_x}\ar[d]_{h_\lambda}& S_x\ar[d]^{h_\lambda}\cr S_y\ar[r]^{J_y}& S_y}$$
\end{enumerate}

Going back and forth from the infinitesimal to the integrated properies is not difficult using the formula

$$(\nabla_X \Psi)(x)=\lim_{t\rightarrow 0}{h_\lambda^{-1}(\Psi(\lambda(t))-\Psi(x)\over t}$$

where $\lambda : [0;1]\rightarrow M$ is a curve such that $\lambda(0)=x$ and $({d\over dt}\lambda)(0)=X$, and the similar formula for $\nabla^{LC}$. The integrated versions continue to make sense in the discrete context, as we will see.

%
%
%


%
%

\begin{remark}
It is immediate that the holonomy operators of a Clifford (resp. spin) connection    belong to the Clifford (resp. spin) group. A Clifford connection which is metric for a $c$-compatible spinor metric would have its holonomies in the spin-c group and would deserve to be called a spin-c connection. Such objects play a role in the theory of Seiberg-Witten invariants and are indeed called that way. Clifford connections were defined in \cite{bgv}.
\end{remark}
\smallbreak

Two connections on $S$ differ by an $End(S)$-valued 1-form, and it is easy to see that this 1-form must be scalar-valued in the case of Clifford connections. If the connections are also metric, the 1-form has values in $i\RR$, and if they commute with $J$ it has values in $\RR$, so we obtain the uniqueness of the spin connection. Moreover the spin connection always exists (see theorem 9.8 in \cite{gracia} for the Riemannian case, \cite{nadir} for the general case).

Every Clifford connection gives rise to a Dirac operator which will have the local form $D=-i\sum_\mu \gamma(dx^\mu)\nabla_{\partial_\mu}$, where $\gamma : T^* M\rightarrow End(S)$ is the composition of   the musical isomorphism defined by $g$,  and the representation $\rho$ (hence $\gamma(dx^\mu)\psi=\rho((dx^\mu)^\sharp)\psi$).

In noncommutative geometry we will be left with (an abstract version of) the sections of the spinor bundle, Krein product, Dirac operator and spin structure, with no direct hold on the connection, Clifford algebra and vector fields. The following facts are then crucial:

\begin{enumerate}
\item The  the Dirac operator commutes with the spin structure iff the latter commutes with the Clifford connection and anti-commutes with vector fields.
\item The Dirac operator is essentially self-adjoint with respect to the Krein product $(.,.)_H$ iff the spinor metric $H$ is $c$-compatible and the Clifford connection is metric  with respect to it.
\end{enumerate}
 
We prove these claims in appendix \ref{append}.

We can conclude from this that the Dirac operator of a Clifford connection $\nabla$ is essentially self-adjoint and commutes with $J$ iff $\nabla$ is the spin connection. We will recover a discrete analog of this result in section \ref{sec54} of \cite{part2}.
 
\subsection{The canonical spectral ``triple'' of a semi-Riemannian manifold}

In this section we consider a space and time orientable spin manifold $M$, with a given spinor bundle $S$ equipped with  $c$-compatible spinor metric $H$ and representation $\rho$, charge conjugation $J$, and the canonical Dirac operator $D$ corresponding to the spin connection $\nabla$.

\subsubsection{The Krein space of spinor fields}

We let $\Gamma^\infty_c(S)$ denotes the space of smooth sections of $S$ with compact support. On this space we have already defined the hermitian form $(.,.)_H$. Since $(M,g)$ is space orientable we can consider a space-orientation\footnote{Once again, in the antilorentzian convention this is a time-orientation !} form $\beta$, and the corresponding smooth field of Clifford group elements $b=\theta(\beta)$. In fact we can suppose without loss of generality that $b_x$ lies in the Pin group for all $x$, since here $c(b_x)=b_x$ and we have $b_x^2\in \RR$. Hence we can replace $b$ with $b/\sqrt{|b^2|}$ (see proposition \ref{admreal} for details).   We set $B=\rho\circ b$, and we recall that we can define a definite spinor metric   $(.,.)_b$ by premultiplying spinors by $B^{-1}$ or $iB^{-1}$ according to the negativity index of the metric $g$ (see proposition \ref{gprod}). For instance in case $B=B^\times$

\be 
(\Psi,\Phi)_b=\int_MH_x(\Psi(x),B_x^{-1}\Phi(x))\sqrt{-g}dx\nonumber
\ee

is a definite hermitian form on $\Gamma^\infty_c(S)$. Since $M$ is connected we can suppose that it is positive definite up to an overall change of sign when choosing $\beta$. We call $K$ the Hilbert space obtained by completing $\Gamma^\infty_c(S)$ with respect to this scalar product. It turns out that this completion does not depend on the choice of $\beta$ (see \cite{baum} section 3.3.1 or \cite{nadir}). Moreover, as we have already noticed at the end of section \ref{sec23}, since $Ad_b\circ c$ is an Euclidean real structure, we have $b_x^2=1$ if $b_x^\times=b_x$ and $b_x^2=-1$ if $b_x^\times=-b_x$. In either case $B^2=1$. We then obtain canonically  a Krein space $K$ with fundamental symmetry $B$.

The canonical spectral ``triple'' of a Riemannian spin manifold of even dimension is the following bunch of objects : the pre-$C^*$-algebra ${\cal A}={\cal C}^\infty(M)$, Hilbert space of $L^2$-spinor fields ${\cal H}$, representation $\pi$ of ${\cal A}$ on ${\cal H}$ by pointwise multiplication, canonical Dirac operator $D$, charge conjugation operator $J$ and chirality operator $\chi$. We see that in the semi-Riemannian space and time orientable context we must replace ${\cal H}$ with the Krein space ${\cal K}$ described above, but without fixing a particular fundamental symmetry: doing so would be equivalent to fix a particular space orientation. We will see later that the $C^*$-structure on ${\cal A}$ has to be dropped too. But for the moment we turn to the question of the characterization of the signature in the Lorentzian and antilorentzian cases.

\subsubsection{The antilorentzian and Lorentzian cases}
  
Given what we have done just above, characterizing the Lorentzian/antilorentzian signature of a manifold metric purely in terms of the data available in noncommutative geometry is easy.   Beware that we now, and for the rest of the paper, revert to the usual convetion of calling ``timelike'' the vectors such that $g(v,v)>0$ in the antilorentzian case. Recall that if $\omega$ is a 1-form we denote by $\gamma(\omega)$ the Clifford multiplication by $\omega$, that is $\gamma(\omega)=\rho(\gamma^\sharp)$.

\begin{theorem}\label{carac1}
Let $(M,g)$ be a semi-Riemannian space and time orientable spin manifold of even dimension, with given   $c$-compatible spinor metric $H$, spin structure $J$ and chirality operator $\chi$. Let $K$ be the Krein space of spinor fields equipped with the Krein product $(.,.)_H$. Then $(M,g)$ is
\begin{enumerate}
\item   antilorentzian iff there exists a   never vanishing 1-form $\beta$ such that  $J\beta J^{-1}=-\beta$ and $(.,\gamma(\beta)^{-1} .)_H$ is  positive definite.
\item  Lorentzian iff there exists a never vanishing 1-form $\beta$ such that  $J\beta J^{-1}=-\beta$ and $(., \gamma(\beta)^{-1}\chi . )_H$ is positive definite.
\end{enumerate}
\end{theorem}

\subsubsection{Causality conditions}\label{causcond}

A very specific feature of the Lorentzian/antilorentzian signature is the existence of a local causal structure. Two approaches which aim at using precisely this feature in noncommutative geometry are reviewed in \cite{besreview}. Though we will not follow this strategy, it is natural to interpret one of the many \emph{global} causality conditions one can put on a spacetime in the terms of theorem \ref{carac1}, namely \emph{stable causality}.

Stable causality  is the fourth stronger causality condition distinguished in \cite{mingsan}, with global hyperbolicity on top. Intuitively a spacetime $(M,g)$ is \emph{stably causal} if $(M,g_\epsilon)$ is causal for every small enough perturbation $g_\epsilon$ of the metric. For the rigorous definition see \cite{mingsan}. Any  simply connected $2$-dimensional spacetime satisfies this property. The following theorem gives a caracterization of stable causality in terms of some classes of functions.

\begin{theorem}\label{carac2} (\cite{san2}) Let $(M,g)$ be a spacetime. Then the following are equivalent
\begin{enumerate}
\item There exists a time function.
\item There exists a temporal function.
\item $(M,g)$ is stably causal.
\end{enumerate}
\end{theorem}



Recall that a \emph{time function} is a continuous strictly increasing function, and a \emph{temporal function} is a smooth function on a spacetime with a past-directed timelike gradient (\cite{mingsan}, def 3.48). 

\begin{remark}
The consequences of this theorem for both noncommutative geometry and quantum gravity might prove far-reaching. First, since causality is defined through (noncommutative) time functions in noncommutative geometry (\cite{besreview}), the equivalence of the first two points is important since it shows that one looses nothing when using the pre-$C^*$-algebra ${\cal C}^\infty(M)$ instead of its $C^*$-completion ${\cal C}(M)$. Moreover, the resilience of  a stably causal metric under small   fluctuations could mean that this condition alone stands a chance of defining a notion a causality in the quantum regime.
\end{remark}
\smallbreak

Putting together theorems \ref{carac1} and \ref{carac2}, one sees immediately that stably causal spacetimes are those on which the $1$-form $\beta$ can be chosen to be exact.

\subsection{Wick rotation of the gamma matrices and Dirac operator}\label{WRglobal}

We take the same hypotheses and notations as in the previous subsection,  except that it will be more convenient here to work with the ungraded charge conjugation operator ${\cal C}$. Also it will be important  to consider a  normalized  space-orientation form $\beta$, so that $BB^\times=1$ and $B^2=\pm 1$.

We know then that $\sigma=Ad_b\circ c$ is a Euclidean real structure   on the Clifford bundle $\CC l(M)$, and we have  seen how $\sigma$  can induce a Wick rotation of the metric. Namely, we have a new metric $g_\sigma(.,.)=g(\sigma(.),.)$, and an isometry $u_\sigma={1+i\over 2}\id+{1-i\over 2}\sigma$ from $(TM,g_\sigma)$,   to the subbundle ${\cal V}_\sigma$ of $\CC l(M)$ equipped with the quadratic form $Q(w)=w^2$. This isometry extends to an isomorphism $\tilde u_\sigma : Cl(M,g_\sigma)\rightarrow Cl_\sigma=\{(x,a)\in \CC l(M)|\sigma(a)=a\}\subset \CC l(M)$.  What we would like to do is to use this isometry to ``Wick rotate'' the Dirac operator directly, short-circuiting the metric, since we do not have a direct access to the metric in noncommutative geometry. The strategy, which, as far as we know, was first used in section 3.4 of \cite{thesevdd}, is to Wick rotate the gamma matrices appearing in the Dirac operator.

We note first that $\rho\circ \tilde u_\sigma$ is a representation of $Cl(M,g_\sigma)$.  Moreover ${\cal C}_\sigma:=B{\cal C}={\cal C}B$  is an ungraded charge conjugation operator corresponding to $\sigma$ (see proposition \ref{generalWR}).

Let us use $\sigma$ to decompose the tangent bundle $TM$ into a $g$-orthogonal sum $TM_+\oplus TM_-$, corresponding at each point to the decomposition given in lemma \ref{decadix} (hence $v$ is a section of $TM_+$ iff $\sigma(v)=v$). Then if $(e_\mu)$ is a local $g$-pseudo-orthonormal frame such that $e_\mu(x)$ lies either in $T_xM_+$ or in $T_xM_-$, then $(e_\mu)$ will also be a $g_\sigma$-pseudo-orthonormal frame.

We can then define gamma matrices with respect to $\rho$ and $\rho\circ\tilde u_\sigma$ :

\begin{itemize}
\item $\gamma_\mu:=\rho(e_\mu)$, which satisfy $\{\gamma_\mu,\gamma_\nu\}=2g(e_\mu,e_\nu)$,
\item $\gamma_\mu^\sigma:=\rho(u_\sigma(e_\mu))$ which satisfy  $\{\gamma_\mu^\sigma,\gamma_\nu^\sigma\}=2g_\sigma(e_\mu,e_\nu)$
\end{itemize}


Since $\rho(\sigma(a))={\cal C}_\sigma \rho(a){\cal C}_\sigma^{-1}$, we have

$$\gamma_\mu^\sigma={1+i\over 2}\gamma_\mu+{1-i\over 2}{\cal C}_\sigma \gamma_\mu{\cal C}_\sigma^{-1}={1+i\over 2}\gamma_\mu+{1-i\over 2}B \gamma_\mu B^{-1}$$

Let $\times_\sigma$ be the adjunction in $End(S)$ for the $\sigma$-compatible Krein product given by proposition \ref{gprod}. We use the same symbol to denote the adjunction in $End(K)$ with respect to $(.,.)_b$.  Then for any operator $A$ on $S$ or $End(K)$  we have $A^{\times_\sigma}=BA^\times B^{-1}$.   Hence we can also write

$$\gamma_\mu^\sigma={1+i\over 2}\gamma_\mu+{1-i\over 2}\gamma_\mu^{\times_\sigma}$$

We must also consider the gamma matrices $\gamma^\mu=\rho((e_\mu^*)^\sharp)=g(e_\mu,e_\mu)\gamma_\mu$, and similarly $\gamma_\sigma^\mu=g_\sigma(e_\mu,e_\mu)\gamma_\mu^\sigma$. We have

$$\gamma_\sigma^\mu=g_\sigma(e_\mu,e_\mu) \big({1+i\over 2}\gamma_\mu+{1-i\over 2} \gamma_\mu^{\times_\sigma}\big)$$

Observe that $g_\sigma(e_\mu,e_\mu)=g(\sigma(e_\mu),e_\mu))=\kappa(\mu)g(e_\mu,e_\mu)$, where $\kappa(\mu)=\pm 1$ according to whether $e_\mu\in T_xM^\pm$. Since $\rho(\sigma(e_\mu))=\gamma_\mu^{\times_\sigma}$, we have 

$$\gamma_\sigma^\mu=g(e_\mu,e_\mu) \big({1+i\over 2}\gamma_\mu^{\times_\sigma}+{1-i\over 2} \gamma_\mu\big)$$

hence

\be
\gamma_\sigma^\mu={1+i\over 2}(\gamma^\mu)^{\times_\sigma}+{1-i\over 2} \gamma^\mu\label{formule11}
\ee



Now we want to replace the matrices $\gamma^\mu$ appearing in $D$ by $\gamma^\mu_\sigma$. In view of (\ref{formule11}), we define

\be
D_\sigma={1+i\over 2}D^{\times_\sigma}+{1-i\over 2}D\label{wickdirac}
\ee

which can also be written

\be
D_\sigma={1+i\over 2}BDB^{-1}+{1-i\over 2}D\label{wd2}
\ee

and from ${\cal C}D=-D{\cal C}$ and ${\cal C}_\sigma=B{\cal C}$ we also obtain

\be
D_\sigma=-{1+i\over 2}{\cal C}_\sigma D{\cal C}_\sigma^{-1}+{1-i\over 2}D\label{wd3}
\ee

It is clear on  expression  (\ref{wickdirac}) that $D_\sigma^{\times_\sigma}=D_\sigma$. From (\ref{wd3}) we get

\bea
\{D_\sigma,{\cal C}_\sigma\}&=&-{1+i\over 2}{\cal C}_\sigma D{\cal C}_\sigma^{-1}{\cal C}_\sigma+{1-i\over 2}D{\cal C}_\sigma-{\cal C}_\sigma{1+i\over 2}{\cal C}_\sigma D{\cal C}_\sigma^{-1}+{\cal C}_\sigma{1-i\over 2}D\cr
&=&-{1+i\over 2}{\cal C}_\sigma D +{1-i\over 2}D{\cal C}_\sigma-{1-i\over 2}{\cal C}_\sigma^2 D{\cal C}_\sigma^{-1}+{1+i\over 2}{\cal C}_\sigma D\cr
&=&0\label{eq16}
\eea

(Hence $[D_\sigma,J_\sigma]=0$.)

Note that in this computation we used the fact that $B^2$ is   constant. 






Finally, if we compute the action of $D_\sigma$ on a spinor field from (\ref{wd2}), we obtain 

\bea
\left({1+i\over 2}BDB^{-1}+{1-i\over 2}D\right)\Psi&=&{1+i\over 2}iB\gamma^\mu\nabla_{e_\mu}(B^{-1}\Psi)+i{1+i\over 2}\gamma^\mu\nabla_{e_\mu}\Psi\cr
&=&i\gamma^\mu_\sigma\nabla_{e_\mu}\Psi+{1+i\over 2}iB\gamma^\mu(\nabla_{e_\mu}^{LC}B^{-1})\Psi\cr
&=&i\gamma^\mu_\sigma\partial_\mu\Psi+\ldots\label{calcul}
\eea
 
In the last line of this expression we have stressed that even though the Wick rotated Dirac operator $D_\sigma$ is not, except for very special metrics, the Dirac operator of the Wick rotated metric $g_\sigma$, it has the same principal symbol. In particular, the metric on $M$ defined by Connes' distance formula from $D_\sigma$ is $g_\sigma$.

The conclusion of this discussion is that $D_\sigma$ has all the algebraic properties required to form a spectral triple, together with the Hilbert space $(K,(.,.)_b)$, charge conjugation $J_\sigma$, chirality $-\chi$, and the representation $\pi$ of ${\cal C}^\infty(M)$ by pointwise multiplication. (We have been very uncaring about the analytical properties. For more on this question see \cite{vddpr}). 

Hence if we set $S=({\cal C}(M),K,(.,.),\pi,D,J,\chi)$ and $S_\beta=({\cal C}(M),K,(.,.)_b,\pi,$ $D_\sigma,J_\sigma,-\chi)$ the map $W_\beta : S\mapsto S_\beta$ is defined   without reference to the Clifford bundle  \emph{as long} as we know how to let $\beta$ act on $K$. In the Lorentzian and antilorentzian cases this action can indeed be defined thanks to $\pi$ and $D$ only. The target of $W_\beta$ is a spectral triple, and the nature of the source is to be defined in \cite{part2}.

 If we want to go in the other direction, we must provide a distinguished form $\beta$ (or equivalently $b$, $c=Ad_{b}\circ \sigma$ or $J=B^{-1}J_\sigma$). From 

$$u_\sigma^{-1}={1-i\over 2}\id+{1+i\over 2}c$$


  we guess that

\be 
 D={1-i\over 2}D_\sigma^\times+{1+i\over 2}D_\sigma={1-i\over 2}B^{-1}D_\sigma B+{1+i\over 2}D_\sigma\label{wickdirac2}
\ee

This is is readily checked (with the first formula we need to use the fact that $(D^{\times_\sigma})^\times=(D^\times)^{\times_\sigma}$ which is true because $BB^\times$ is constant, and with the second formula we need to use $BDB^{-1}=B^{-1}DB$ which is true because $B^2$ is constant). 












Example : We start from $M=\RR^{2k}$ with the constant Euclidean metric. The spinor bundle is trivial and identified with $M\times \CC^{2^k}$. We let $e_\mu=\partial_\mu$. To Wick rotate to antilorentzian signature we take $g=e_0$ (for instance). Then $V_\sigma$ is the vector subspace generated by $e_0$ and $ie_a$, $a=1..2k$, and $G=\gamma_0$ is constant. Then (\ref{wickdirac}) gives

\bea
D_\sigma&=&{1+i\over 2}i(\gamma_0\partial_0-\sum_a\gamma_a\partial_a)+{1-i\over 2}i(\gamma_0\partial_0+\sum_a\gamma_a\partial_a)\cr
&=&i(\gamma_0\partial_0-\sum_a (i\gamma_a)\partial_a)\cr
&=&i(\tilde\gamma_0\partial_0-\sum_a(\tilde\gamma_a\partial_a)),\mbox{with } \tilde\gamma_\mu\mbox{ antilorentzian gamma matrices}\cr
&=&i\sum_a\tilde\gamma^a\partial_a,\mbox{ raising the indices with the antilorentzian metric}\cr
&=&D_{antilor}
\eea

In this simple example $D_\sigma$ \emph{is} the canonical Dirac operator of $g_\sigma$, since the Christoffel symbols vanish and $B$ is constant. 

\appendix

\section{Krein products on algebraic spinors}\label{kreinalgspin}

\begin{propo}\label{propo8} Let $S_e$ be an algebraic spinor module.
\begin{enumerate}
\item If $ee^{\times_\sigma}=0$ then $(.,.)_\sigma$ is zero on $S_e$, hence is not a Krein product.
\item If $ee^{\times_\sigma}\not=0$ then : 
\parbox[t]{8cm}{\begin{enumerate}
\item There exists a unique primitive idempotent $f\in \CC l(V)$ such that $f^{\times_\sigma}=f$ and $S_e=S_f$.
\item The left representation $L$ of $\CC l(V)$ on $(S_e,(.,.)_\sigma)$ is $\sigma$-compatible.
\end{enumerate}}
\end{enumerate}
\end{propo}

\begin{demo}
Suppose $ee^{\times_\sigma}=0$. Then for all $a,b\in\CC l(V)$, $(ae,be)_\sigma=\tau_n(e^{\times_\sigma} a^{\times_\sigma} be)=\tau_n(a^{\times_\sigma} bee^{\times_\sigma})=0$.

Now suppose $ee^{\times_\sigma}\not=0$. Let $(.,.)_a$ be a $\sigma$-compatible Krein spinor product on $S_e$, which exists by proposition (\ref{rob1}).  

Since $e$ is a primitive idempotent, $p:=L_e\in End(S_e)$ is a rank one projection. Note that  $p(e)=e^2=e$, hence $e\in \im(p)$, so that $\im(p)=\CC e$. Observe also that $(\ker p)^\perp=\im(p^\times)$ (the orthogonal and adjoint are relative to $(,)_a$) is one-dimensional\footnote{Proof of these relations : see appendix.}. Let us call $v$ a non-zero vector in this space.

If $v$ were isotropic, then for all $\psi,\eta\in S_e$ one would have $(p^\times(\psi),p^\times(\eta))_a=(\psi,pp^\times(\eta))_a=0$. By the non-degeneracy of $(.,.)_a$ this would imply $pp^\times=0$ which is excluded since $ee^{\times_\sigma}\not=0$.

Then we can suppose that $(v,v)=\epsilon=\pm 1$. Let us call $q\in S_e$ the orthogonal projection on $\CC v=(\ker p)^\perp$. Explicitly one has $q(\psi)=\epsilon(v,\psi)v$, and it is easy to check that $q^2=q=q^\times$. Moreover, $\ker q=v^\perp=\ker p$. Since $v$ is not isotropic, $v\notin \ker(q)$, hence  $S_e=\ker(q)\oplus \CC v$. If $\psi\in \ker(q)$, then $\psi\in\ker p$, hence $pq(\psi)=p(\psi)=0$. Moreover $pq(v)=p(v)$. Thus $p=pq$. Similarly it is easy to check that $qp=q$ using $S_e=\ker p\oplus \CC e$.

Let us call $f=L^{-1}(q)\in \CC l(V)$. Using $L$ one sees that $f^2=f$, $f^{\times_\sigma}=f$, $e=ef$, $f=fe$. This shows a). The uniqueness is easy to prove.

Let us show b). Since we know that $(a\psi,\phi)_\sigma=(\psi,a^{\times_\sigma}\phi)_\sigma$ for all $a\in \CC l(V)$ and $\phi,\psi\in S_e$ by (\ref{kreinmultadj}), then by lemma \ref{unikrein} there exists $\lambda\in\RR$ such that $(.,.)_\sigma=\lambda(.,.)_a$ on $S_e$. There thus suffices to show that $\lambda\not=0$.

Let us consider $f$ as an element of $S_e$. Then $(f,f)_\sigma=\tau_n(f^{\times_\sigma} f)=\tau_n(f)$. Now $f$ being a non zero projection, its normalized trace is not zero. On the other hand we have $(f,f)_\sigma=\lambda (f,f)_a$. This shows that $\lambda\not=0$, and that $(.,.)_\sigma$ is $\sigma$-compatible.
\end{demo}

Let  $S_f$ be an algebraic spinor module with $f$  a primitive idempotent such that $f^{\times_\sigma}=f$. Then for each $x=af$ and $y=bf$ one has $x^{\times_\sigma}y=fa^{\times_\sigma}bf\in f\CC l(V)f$. But remember that in any   finite-dimensional algebra $A$ over $k$, if $f$ is a primitive idempotent then $fAf$ is a division algebra over $k$. Hence we see that the range of the map $(x,y)\mapsto x^{\times_\sigma}y$ defined on $S_e\times S_e$ is a division algebra over $\RR$, hence is isomorphic to $\RR$, $\CC$ or ${\mathbb H}$. Such a map is called\footnote{We thank Christian Brouder for drawing our attention to this reference.} a \emph{spinor product} in \cite{lou}. Thus we see that the $\sigma$-compatible Krein product on algebraic spinor modules are just the composition of Lounesto's spinor product with the normalized trace. It should certainly be possible to relate the properties of these two kinds of products.


\begin{remark} Given a $\sigma$-compatible representation $\rho$ in a Krein space $K$, we see that $\rho(S_e)$ is the set of operators $a$ such that $\im(a^\times)\subset L$ where $L$ is a line. The set of minimal left ideals $S_e$ has thus the structure of $P(\CC^{2^{n/2}})$. Moreover, the ideals on which $(.,.)_\sigma$ vanishes are exactly those whose corresponding line is isotropic. They form a submanifold of the projective space which is the image by the quotient map of the isotropic cone of the Krein product.
\end{remark}

\section{$C^*$-algebra structures on $\CC l(V)$}\label{cstarstruccliff}

Let $\sigma$ be an Euclidean real structure.  Let $K$ be an irreducible spinor module and $\bra .,.\ket_\sigma$ a $\sigma$-compatible scalar product on $K$. Then let 
$$\|a\|_\sigma^\rho:=\sup_{\psi\not=0}{\bra a\psi,a\psi\ket_\sigma\over \bra \psi,\psi\ket_\sigma}$$
Obviously, $(\CC l(V),\times_\sigma,\|.\|_\sigma^\rho)$ is a $C^*$-algebra isomorphic to $M_{2^{n/2}}(\CC)$ with the usual involution $*$ and operator norm $\|.\|$ associated with the canonical scalar product. Hence if we equip   $\CC l(V)$ with  $\sigma$ we obtain  an algebra with two commuting involutions, $\times$ and $\times_\sigma$, one of them being a $C^*$-involution. This kind of structure is called a \emph{Krein $C^*$-algebra} in \cite{kawa} (see also \cite{berto}). Thus, defining a Wick rotation of the metric to Euclidean signature exactly corresponds to giving $\CC l(V)$ a Krein $C^*$-algebra structure. Note that we used a particular representation to make this remark (in example 4.10 of \cite{kawa}, where this point is also emphasized, the natural representation of the Clifford algebra on the exterior algebra is used). However if we use the $\sigma$-product there is no need to leave the Clifford algebra.

Indeed, notice that the $\sigma$-product $(a,b)_\sigma=\tau_n(a^{\times_\sigma}b)$ is transported by $\rho$ to the Hilbert-Schmidt scalar product on $M_{2^{n/2}}(\CC)$. In particular one has $\|a\|_\sigma^2=\tr(\rho(a)^*\rho(a))$. Now write $\|A\|_{HS}:=[\tr(A^*A)]^{1/2}$ for the Hilbert-Schmidt norm of a matrix $A$ and define $\|A\|_{\infty,HS}=\sup\{\|AB\|_{HS}\ |\ \|B\|_{HS}=1\}$. Then it is easy to see that $\|A\|_{\infty,HS}$ is in fact equal to $\|A\|$. Hence the $C^*$-norm on $\CC l(V)$ corresponding to the involution $\times_\sigma$ can be defined without using the representation $\rho$ by 

\be
\|a\|_{\infty,\sigma}:=\sup\{\|ab\|_\sigma\ |\ \|b\|_\sigma=1\}\nonumber
\ee

Conversely, if $*$ is a $C^*$-involution which commutes with $\times$ and stabilizes $V^\CC$ then it is of the form $\times_\sigma$ for some Euclidean real structure $\sigma$. Moreover if $\sigma'$ is another Euclidean real structure, then $(\CC l(V),\times_\sigma,\|.\|_{\infty,\sigma})$ and $(\CC l(V),\times_{\sigma'},\|.\|_{\infty,\sigma'})$ are isomorphic. In fact this isomorphism is canonical and is just $Ad_{b{b'}^{-1}}$.

\section{Proofs of the claims in subsection \ref{sec33}}\label{append}
Let $(M,g)$ be a semi-riemannian spin manifold, with spinor bundle $S$ and spin structure $J$. Let $\nabla$ be a Clifford connection and $D$ be the Dirac operator $D$ associated to it. We write $\nabla_\mu:=\nabla_{\partial\over\partial\mu}$. 

\begin{propo} \label{RealStructureThm}
The Dirac operator and spin structure commute if and only if $\{J, X\} = [J,\nabla_X] = 0$ for all $X \in \Gamma(TM)$.
\end{propo}

\begin{demo}
The fact that $\{J,X\} = [J,\nabla_X] = 0$ are sufficient conditions is a straightforward calculation.

Let us prove that they are necessary. Let $f$ be a smooth real-valued function. Since $[J,  {D}]=[J,f]=0$, we deduce that $J$ commutes with $[ {D}, f] = -i \gamma(df)$. We infer easily that $J$ anti-commutes with any real vector field or differential form. Next we consider the differential operator $\nabla'= J^{-1} \nabla  J$. This operator is a Clifford connection. Indeed, for any $X,Y \in \Gamma(TM)$, $\psi \in \Gamma(S)$:

\bea
\nabla_X'(Y\cdot \psi) &=& J^{-1} \nabla_X  J(Y\cdot \psi) \cr
&=& -J^{-1} \nabla_X  ( Y \cdot J\psi) \cr
&=& -J^{-1} (\nabla_X^{LC} Y)\cdot J\psi + Y\cdot \nabla_X  (J\psi) \cr
 &=&  (\nabla_X^{LC} Y)\cdot \psi + Y\cdot\nabla_X' \psi
\eea

There thus exists a complex-valued one-form $\omega$ such that: $\nabla_X' - \nabla_X= \omega(X)$. This can also be written: $\omega(X) = J^{-1} [\nabla_X, J]$. We have:
\bea
0 &=& [J,  {D}] \cr
&=& [J, -i\sum_\mu\gamma(dx^\mu) \nabla_\mu] \cr
&=& -i\sum_\mu\gamma(dx^\mu)[J, \nabla_\mu] \cr
 &=& -i J \sum_\mu\gamma(dx^\mu) \omega_\mu \nonumber
\eea

This implies that $\omega=0$, and thus that $[\nabla_X, J]=0$.
\end{demo}

Let $H$ be a spinor metric and $(.,.)_H$ be the associated product on spinor fields defined by (\ref{kreinprod}).

\begin{propo}
The Dirac operator is essentially self-adjoint if and only if $X^\times =X$ for all $X \in \Gamma(TM)$, and the Clifford connection $\nabla$ is metric for $H$.
\end{propo}

\begin{remark} Our spinor bundle $S$ thus has to be a Dirac bundle (see \cite{LM}, p 114).
\end{remark}

\begin{demo}
Let us assume that the Dirac operator is self-adjoint. Let $f$ be a real-valued smooth function. Then $f$ is self-adjoint. This implies that $[ {D}, f] = -i \gamma(df)$ is anti-self-adjoint. We easily conclude from this that all real vector fields and all real-valued 1-forms are self-adjoint. In particular, the $\gamma(dx^\mu)$ are self-adjoint. Let $\psi, \phi$ be   spinor fields with compact support. We have:
\bea
0 &=& (\psi,  {D} \phi)_H - ( {D} \psi, \phi)_H \cr
&=& \int \sqrt{|g|}\sum_\mu \left[H(-i\gamma(dx^\mu) \nabla_\mu\psi, \phi) - H(\psi, -i\gamma(dx^\mu) \nabla_\mu \phi)\right]dx \cr
0 &=& i\int \sqrt{|g|}\sum_\mu \left[H(\nabla_\mu \psi, \gamma(dx^\mu) \phi) + H(\psi, \gamma(dx^\mu) \nabla_\mu \phi)\right]dx \nonumber
\eea
We also have that: 
\bea
\left[\nabla_\mu, \gamma(dx^\mu)\right] &=& \gamma(\nabla_\mu^{LC}(dx^\mu)), \cr
&=& -\sum_\alpha\Gamma^\mu_{\mu\alpha} \gamma(dx^\alpha)\nonumber
\eea
Hence
\bea
\sum_\mu\left[\nabla_\mu, \gamma(dx^\mu)\right]&=& -\sum_{\mu,\alpha}\Gamma^\mu_{\mu\alpha} \gamma(dx^\alpha)\cr
&=& -\sum_\alpha\frac{\partial_\alpha \sqrt{|g|}}{\sqrt{|g|}}\gamma(dx^\alpha)\nonumber
\eea
from which we infer that $\sum_\mu\gamma(dx^\mu) \nabla_\mu  = \sum_\mu\nabla_\mu  \gamma(dx^\mu) +\sum_\mu \gamma(dx^\mu) (\partial_\mu \sqrt{|g|})/ \sqrt{|g|}$. Substituting in the integral above gives us:
$$
\int dx\sum_\mu\left[ \sqrt{|g|}  H(\nabla_\mu  \psi, \gamma(dx^\mu) \phi) +  \sqrt{|g|} H(\psi, \nabla_\mu (\gamma(dx^\mu) \phi))  + (\partial_\mu \sqrt{|g|}) H(\psi, \gamma(dx^\mu) \phi)  \right] = 0
$$
Finally, an integration by part yields:
$$
\int dx \sqrt{|g|}\sum_\mu \left[ H(\nabla_\mu \psi, \gamma(dx^\mu) \phi) + H(\psi, \nabla_\mu (\gamma(dx^\mu) \phi)) - \partial_\mu H(\psi, \gamma(dx^\mu) \phi)  \right]dx = 0
$$
for all $\psi, \phi \in \Gamma_c(S)$. Now, the expression between brackets can be proven to be ${\mathcal C}^\infty (M, \mathbb{C})$-linear in $\phi$ (and anti-linear in $\psi$ as well). Indeed, $H(\nabla_\mu  \psi, \gamma(dx^\mu) \phi)$ is clearly linear in $\phi$. let $f \in {\mathcal C}^\infty (M, \mathbb{C})$. We replace $\phi$ by $f \phi$ in the two remaining terms:

\bea
H(\psi, \nabla_\mu  (\gamma(dx^\mu) f\phi)) - \partial_\mu H(\psi, \gamma(dx^\mu) f\phi) &=& H(\psi, \nabla_\mu  [f (\gamma(dx^\mu) \phi)]) - \partial_\mu [f H(\psi, \gamma(dx^\mu) \phi)] \cr
&=& H(\psi, (\partial_\mu f) \gamma(dx^\mu) \phi) + f \nabla_\mu  (\gamma(dx^\mu) \phi)) \cr
&&- (\partial_\mu f) H(\psi, \gamma(dx^\mu) \phi) - f \partial_\mu H(\psi, \gamma(dx^\mu) \phi) \cr
&=& f [H(\psi, \nabla_\mu (\gamma(dx^\mu) \phi)) - \partial_\mu H(\psi, \gamma(dx^\mu) \phi)]\nonumber
\eea

Thus, for all $f \in {\mathcal C}^\infty (M, \mathbb{C})$ and $\psi, \phi \in \Gamma_c(S)$:
$$
\int dx \sqrt{|g|} f\sum_\mu \left[ H(\nabla_\mu  \psi, \gamma(dx^\mu) \phi) + H(\psi, \nabla_\mu  (\gamma(dx^\mu) \phi)) - \partial_\mu H(\psi, \gamma(dx^\mu) \phi)  \right] = 0
$$
which implies that:
\be
\sum_\mu\left[H(\nabla_\mu \psi, \gamma(dx^\mu) \phi) + H(\psi, \nabla_\mu  (\gamma(dx^\mu) \phi)) - \partial_\mu H(\psi, \gamma(dx^\mu) \phi)\right]=0\label{equationonHdivergence}
\ee

at every $x$. Now suppose that $\nabla^0$ is a local spin connection defined around $x$. We know that $\nabla=\nabla^0+A$ where $A$ is a scalar 1-form, and that $\nabla$ is metric if and only if $A$ has pure imaginary values. Using (\ref{equationonHdivergence}) we find that 

\be
\sum_\mu (A_\mu+\bar A_\mu)H(\psi,\gamma(dx^\mu)\phi)=0\label{eqaz}
\ee

for all spinor fields $\psi$, $\phi$ with small enough compact support containing $x$. Since $H$ is non-degenerate the orthogonal of $\{\gamma(dx^2)\phi,\ldots,\gamma(dx^n)\phi\}$ for $H$ is a subbundle of dimension $2^{n\over 2}-n+1\ge 1$ of the spinor bundle. Since the orthogonal of $\{\gamma(dx^1)\phi,\ldots,\gamma(dx^n)\phi\}$  has codimension 1 in it, we can, at least locally, consider a spinor field $\psi$ such that $H(\psi,\gamma(dx^\mu)\phi)=0$ for $\mu=2,\ldots,n$ and $H(\psi,\gamma(dx^1)\phi)\not=0$. Hence $A_1+\bar A_1=0$, and of course the same can be done for the other indices. The Clifford connection $\nabla$ is thus metric.

The converse can be proven easily following the same steps. It is in fact a standard result of spin geometry. See for example \cite{LM}.
\end{demo}

\end{document}